\documentclass[12pt]{article}

\usepackage{biblatex}
\usepackage{moresize}
\usepackage{graphicx}
\usepackage{caption}
\usepackage{subcaption}
\usepackage{tikz}
\usepackage{amsmath}
\usepackage{amsfonts} 
\usepackage{wrapfig}
\usepackage[english]{babel}
\usepackage{hyperref}
\usepackage{lipsum}
\usepackage{mdframed}
\usepackage{enumitem}
\usepackage{csquotes}

\newtheorem{theorem}{Theorem}[section]

\newtheorem{definition}[theorem]{Definition}

\numberwithin{equation}{section}

\usepackage[margin=2.5cm]{geometry}

\usepackage{newtxtext}

\usepackage{setspace}
\usepackage{parskip}

\title{\vspace{-2cm}Stabilization Approach for Spectral-Volume Methods}
\date{\today}
\author{Lena Schadt}

\begin{document}

\title{Stable Spectral-Volume Methods}
\author{Lena Schadt\footnote{\begin{minipage}[t]{\linewidth}l.schadt@tu-bs.de\\[0.1cm]
Institute for Partial Differential Equations, TU Braunschweig, Germany\end{minipage}}}
    
\date{\today}

\maketitle
\vspace{1cm}

\textbf{Abstract}\\
A novel approach for the stabilization of the Spectral-Volume (SV) method based on Dafermos' entropy rate criterion is presented. The method is an adaption of an already existing approach for the stabilization of the Discontinuous-Galerkin (DG) method. It employs the same estimates for the maximal possible entropy dissipation rate as the DG version. However, a new way to compute the discrete conservative filter had to be derived due to the differences of the underlying schemes. The resulting modified SV scheme even satisfies the entropy inequality. Tests are carried out for Burgers' equation and for the Euler equations of gas dynamics.

\vspace{0.4cm}


\vspace{0.4cm}


\vspace{0.4cm}

\section{Introduction}\label{sec:Introduction}

The partial differential equations we aim to solve numerically are one-dimensional hyperbolic systems of conservation laws
\begin{equation}
    \begin{aligned}\label{eq:conservation law}
        \partial_{t}u\left(x,t\right) +  \partial_{x} f \left(u\left(x,t\right)\right) &= 0, \quad
        &&\left(x,t\right) \in \Omega \times \left[0,\infty \right[ \\
        u\left(x,0\right) &= u_{0}\left(x\right), \quad &&x \in \Omega
    \end{aligned}
\end{equation}
with \( \Omega \subseteq \mathbb{R} \). The function \( u : \Omega \times \left[0, \infty \right[ \longrightarrow    \mathbb{R}^{m} \) is the unknown vector of state variables and \(f : \mathbb{R}^{m}\longrightarrow\mathbb{R}^{m}\) is a given flux function. A function \(u\) is called a \textit{classical solution} of \eqref{eq:conservation law} if \(u\) is a \(C^{1}\)-function that satisfies \eqref{eq:conservation law} pointwise. Unfortunately, it is well known that in many cases classical solutions to conservation laws break down. Consequently, one introduces the concept of weak solutions \cite{godlewski2013numerical,lax1971shock}.

Let \( D = \mathbb{R} \times \left[0, \infty \right[ \), and \( C_{0}^{\infty} \left( D \right) \) denote the space of all smooth functions with compact support. A function \( \varphi\) is termed a testfunction if \( \varphi \in C_{0}^{\infty} \left( D \right) \). For brevity, we denote \(u = u\left(x, t\right) \) and \(\varphi = \varphi\left(x,t\right) \). Multiplying (\ref{eq:conservation law}) by \( \varphi \in C_{0}^{\infty} \left( D \right) \) and integrating over \(D\) yields

\begin{equation}
    \int_{0}^{\infty} \int_{-\infty}^{\infty} \varphi \partial_t u + \varphi \partial_x f\left(u \right) \, dx dt = 0.
\end{equation}
By performing an integration by parts and considering the compact support of \(\varphi\), we obtain

\begin{equation}\label{eq:weak}
    \int_{0}^{\infty} \int_{-\infty}^{\infty}  u \partial _t \varphi + f\left(u \right) \partial_x \varphi \, dx dt = \int_{- \infty}^{\infty} \varphi\left(x, 0 \right) u\left( x, 0 \right)\, dx.
\end{equation}
This formulation is valid even if \(u \in L^{\infty}_{loc} \left( D \right)\) only. A function u is termed a weak solution of \eqref{eq:conservation law} if \eqref{eq:weak} holds for all \(\varphi \in C_{0}^{\infty} \left( D \right)\). We say that \(u\) satisfies \eqref{eq:conservation law} in the weak sense \cite{lax1971shock}.

Unfortunately, weak solutions are not unique [6, example 2.3]. Therefore, we will introduce entropy
conditions as an additional constraint on the solution. These conditions should enable us to pick out the physically relevant weak solutions.

A convex function \(U: \mathbb{R}^{m} \rightarrow \mathbb{R} \) is an entropy, if there exists a function \( F: \mathbb{R}^{m} \rightarrow \mathbb{R}\), called entropy flux, such that
\begin{align}\label{eq:entropy equation}
    U'\left(u\right) f'\left(u\right) = F'\left(u\right)
\end{align}
holds \cite{lax1971shock}. An admissable pair of entropy and entropy flux will be refered to as entropy-entropy flux pair \( \left(U,F\right) \). Note that \(U' = \left(\nabla U \right)^{T}\), \(F' = \left( \nabla F \right)^{T}\) and \(f'\) denotes the jacobian matrix of \(f\).

If a scalar conservation law is considered as with the linear advection equation, every convex function \( U: \mathbb{R} \rightarrow \mathbb{R} \) could be selected as entropy. A corresponding entropy flux would then be given by \cite{lax1971shock}
\begin{equation}
    F\left(u \right) = \int  U'\left(u\right) f'\left(u\right) du.
\end{equation}
Note further that equation (\ref{eq:entropy equation}) is overdetermined if \(m > 2\). Hence, it is possible that no entropy-entropy flux pair exists at all. Fortunately, in the case of the Euler equations of gasdynamics suitable entropy- entropy flux pairs are known \cite{harten1983symmetric, tadmor2003entropy}.

If \(u\) is a classical solution of \eqref{eq:conservation law} it satisfies another conservation law as well. Multiplying \eqref{eq:entropy equation} by \(\partial_x u\) and using the chain rule we obtain

\begin{align}\label{eq:entropy equality}
    \partial_t U\left(u\right) + \partial_x F\left(u\right) = 0
\end{align}
In general, one can show that for an entropy-entropy flux pair \( \left(U, F\right)\)

\begin{align}\label{eq:entropy_inequality}
    \partial_t U\left(u\right) + \partial_x F\left( u \right) \leq 0
\end{align}
holds in the weak sense \cite{lax1971shock}. For \(m = 1\) the additional condition given by equation \eqref{eq:entropy_inequality} enforces the uniqueness of the solution to \eqref{eq:conservation law}  \cite{kruvzkov1970first}, if it holds for a suitable family of entropies.  
For \(m \geq 2\) on the other hand, equation \eqref{eq:entropy_inequality} yields no guaranty for a unique solution \cite{feireisl2020oscillatory}. Therefore, we use a different entropy criterion.

Let \( u\) be a weak solution of \eqref{eq:conservation law} and \(U\) an admissable entropy. The total entropy of \(u\) at time \(t\) is then defined as:

\begin{align}
    E_{u}\left( t \right) = \int U\left(u\left(x, t\right)\right) \, dx.
\end{align}
Dafermos' entropy rate criterion \cite{dafermos1973entropy} states that the total entropy of the relevant weak solution  \(\tilde{u}\) should reduce faster than the total entropy of any other existing weak solution \(u\):

\begin{align}\label{eq:entropy_rate_criterion}
    \forall t > 0: \quad \frac{dE_{\tilde{u}}}{dt} \leq \frac{dE_{u}}{dt}.
\end{align}
See \cite{dafermos2012maximal, dafermos2009variational, feireisl2014maximal} for theoretical examples where Dafermos' criterion is able to reduce the amount of admissible solution. We will use this criterion in chapter \ref{sec:Entropy Stabilization} as basis for the stabilization approach.

\subsection{Short Introduction to the FV Method}
We start by a brief introduction to the Finite-Volume (FV) method \cite{godlewski2013numerical} as a fundamental concept\. The Spectral-Volume method is essentially based on this simpler and better known method.

First, the domain \(\Omega = \left[a, b\right]\) is divided into \(N\) volumes \(V_i\) called cells

\begin{equation}
     \left[a, b \right] = \bigcup_{i=1}^{N} V_{i}, \quad V_{i} = \left[x_{i-1/2}, x_{i+1/2} \right], \quad i = 1,...,N.
\end{equation}
To keep things simple, we assume an equidistant partition and denote the length of one cell by  \( h = x_{i+1/2} - x_{i-1/2},  \, i = 1,...,N\). The center of each cell is given by

\begin{equation}
    x_i = \frac{x_{i-1/2} + x_{i+1/2}}{2}, \quad i = 1,...,N,
\end{equation}
and the cell average of \(u\) in the \(i\)-th volume will be denoted by an upper bar

\begin{equation}
    \bar{u}_i \left(t \right) = \frac{1}{h} \int_{x_{i-1/2}}^{x_{i+1/2}} u\left(x,t\right) dx.
\end{equation}
The FV method can be derived from an integration of equation \eqref{eq:conservation law} over one volume \(V_i\)

\begin{equation}
     \int_{x_{i-1/2}}^{x_{i+1/2}} \partial_t u \left(x, t\right) + \partial_x f \left( u \left(x, t\right) \right) dx = 0.
\end{equation}
If we assume \(u\) to be sufficiently smooth, the fundamental theorem of calculus yields

\begin{equation}\label{eq:FV_intermediate_step}
    h \cdot \frac{d \bar{u}_i}{dt} \left( t \right) + \left( f_{i+1/2} \left( t \right) - f_{i-1/2} \left( t \right) \right) = 0,
\end{equation}
where
\begin{equation}\label{eq:FV_flux}
    f_{i+1/2}\left( t \right) = f \left( u \left(x_{i+1/2}, t \right) \right).
\end{equation}
Equation \eqref{eq:FV_intermediate_step} is equivalent to
\begin{equation}\label{eq:FV_semidiscrete}
    \frac{d \bar{u}_i}{dt} \left( t \right) = \frac{1}{h} \left( f_{i-1/2} \left( t \right) - f_{i+1/2} \left( t \right) \right).
\end{equation}
and specifies the time evolution of the cell averages in all volumes. Intuitively equation \eqref{eq:FV_semidiscrete} states that the amount of a conserved quantity within a volume \(V_i\) only changes by its flux across the boundary. Unfortunately, the equation can't be integrated directly. That is because we have no access to nodal values of \(u\) at volume boundaries and are thus unable to evaluate the flux at these points. Instead, we approximate the flux at volume boundaries by a numerical flux \(f^*\), that depends on approximations of the left-hand and right-hand limits \(u_l\) and \( u_r \) of \(u\) at a given volume boundary. The most simple way to approximate the one-sided limits is to use the cell averages \(u_l = \bar{u}_i\) and \(u_r = \bar{u}_{i+1}\) of the adjacent volumes directly. To get better approximations of \(u_l\) and \(u_r\) often a reconstruction step is added. That implies that for every step in time \(u\) is approximated in each volume by a polynomial of low degree. The coefficients of the polynomial are determined by the cell averages of the surrounding volumes. The one-sided limits \(u_l\) and \(u_r\) are then calculated by an evaluation of the reconstruction polynomials at the volume boundaries.

\subsection{Basic Idea of the SV Method}

As with the FV method, the fundamental concept behind Spectral-Volume methods \cite{wang2002spectral} is to partition the spatial domain into so-called spectral volumes (SV). However, for the SV method each spectral volume is further subdivided into control volumes (CV). Cell-averaged data from these control volumes is used to reconstruct a high-order polynomial approximation within each spectral volume. In contrast to the FV method, the reconstruction polynomial remains continuous throughout the entire spectral volume. Consequently, we can use the analytical flux function within each spectral volume. Only at the SV boundaries do we need to employ numerical fluxes to address the arising Riemann problems.

In this chapter, all vectors of cell averages and nodal values are presented in bold typeface, as are all matrices. Cell averages are denoted by an upper bar, as in \(\bar{u}\). Generally, we will use \(i\) as an index for the SVs and \(j\) as an index for the CVs. Besides, a superscript star will indicate the numerical flux and the numerical entropy flux, as in \(f^{*}\) and \(F^{*}\).

\subsection{Spatial Grid}

\begin{figure}[!h]
    \centering
    \includegraphics[width=\linewidth]{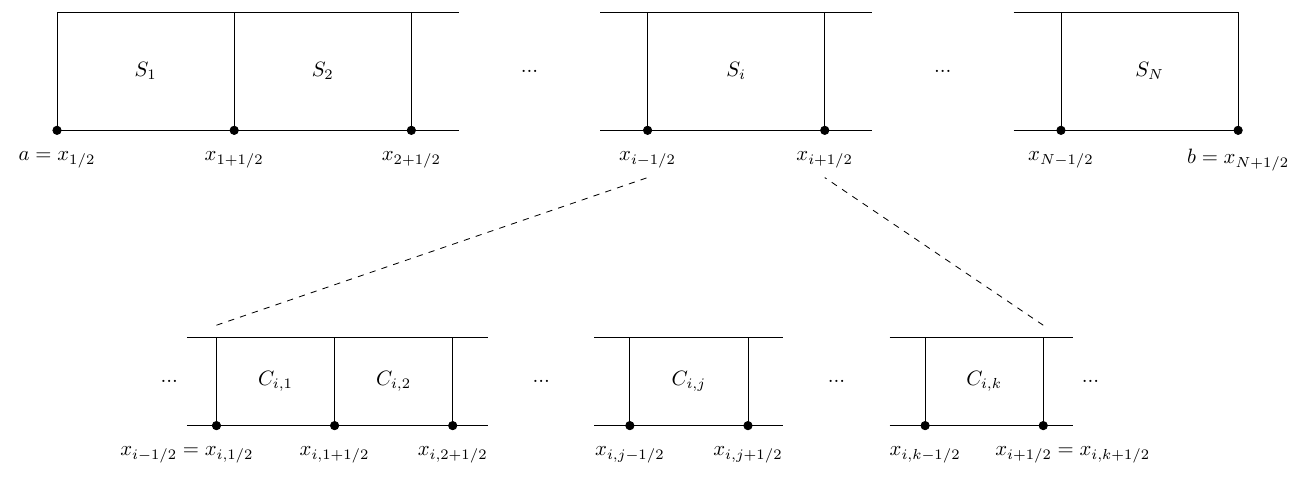}
    \caption{Partition into spectral volumes and further subdivision into control volumes.}
    \label{fig:grid}
\end{figure}

We will now look at the concept of Spectral-Volume methods \cite{wang2002spectral} in more detail. As with the FV method, we opt for a partition \( \{ x_{i+1/2} \} _{i=0}^{N}\) of the domain \( \Omega = \left[a, b\right] \) to divide it into N spectral volumes \(S_{i}\), where \(x_{1/2} = a \) and \(x_{N+1/2} = b\)

\begin{equation}
    \left[a, b \right] = \bigcup_{i=1}^{N} S_{i}, \quad S_{i} = \left[x_{i-1/2}, x_{i+1/2} \right], \quad i = 1,...,N.
\end{equation}
In all numerical tests, we choose an equidistant partition, although this selection is not mandatory. However, as there are no specific requirements for the partition, we will proceed with the special case of equidistant spectral volumes for simplicity. Let \( h = x_{i+1/2} - x_{i-1/2},  \, i = 1,...,N\) represent the length of one spectral volume and
\begin{equation}
    x_i = \frac{x_{i-1/2} + x_{i+1/2}}{2}, \quad i = 1,...,N
\end{equation}
be the center of each SV. The subsequent step involves determining the number of control volumes into which each spectral volume will be subdivided. This decision influences the numerical order of accuracy in space, as we utilize the cell-averages of all control volumes within one spectral volume for the reconstruction. For a desired accuracy of order k, each spectral volume \(S_{i} \) should be divided into \(k\) control volumes using a second partition \( \{ x_{i,j+1/2}\}_{j=0}^{k}\) with \(x_{i,1/2} = x_{i-1/2}\) and \( x_{i,k+1/2} = x_{i+1/2}\). That is because the reconstruction polynomial within each SV will be of order \(k-1\) if \(k\) control volumes are selected. The \(j\)th CV of \(S_{i}\) will be denoted by

 \begin{equation}
    C_{i,j} = \left[ x_{i, j-1/2}, x_{i, j+1/2} \right],  \quad i = 1,...,N, \quad j = 1,...,k
\end{equation}
and \( h_{j} = x_{i,j+1/2} - x_{i,j-1/2}, \, i = 1,...,N\) will be the length of the jth control volume. A visualization of the construction of the grid is depictured in figure \ref{fig:grid}.

In contrast to the partition of \( \left[a, b\right]\) into spectral volumes, we aim to employ a non-equidistant partition to subdivide each spectral volume into control volumes. This choice is motivated by the potential for the reconstruction polynomial to exhibit strong oscillations when equal length control volumes are utilized. This phenomenon is well-known in polynomial interpolation on equidistant grid points and is mitigated by adopting a grid that becomes denser towards the edges. Therefore, in all numerical tests, Gauss-Lobatto points \cite{hesthaven2007nodal, gassner2013skew} are used for the subdivision of the SVs into CVs.

Using a linear transformation, we obtain the partition  \( \{ x_{i,j+1/2}\}_{j=0}^{k}\) from the Gauss-Lobatto points \( \tilde{x}_{j+1/2} \) by
\begin{equation}\label{eq:lintrans}
    \begin{aligned}
        x_{i,j+1/2} = x_i + \frac{h \cdot \tilde{x}_{j+1/2}}{2}, \quad i = 1,...N, \quad  j = 0,...,k,
    \end{aligned}
\end{equation}
where \(h\) was defined as the length of one SV.

\subsection{Polynomial Reconstruction}\label{sec:Polynomial Reconstruction}

Now that we've established the grid for our computations, let's delve into the core of the spectral volume method—polynomial reconstruction.
At the onset of each numerical test, the initial condition \( u\left(x,0\right) = u_{0}\left(x\right) \) needs to be transformed into a cell-averaged initial vector \(\mathbf{\bar{u}}\) with components
\begin{equation}
    \begin{aligned}
        \bar{u}_{i,j} = \frac{1}{h_j} \int_{x_{i,j-1/2}}^{x_{i,j+1/2}} u_0\left(x\right) dx, \quad i = 1,...,N, \quad j = 1,...,k.
    \end{aligned}
\end{equation}
The integral can be approximated numerically, for exammple, using Gaussian quadrature \cite{isaacson2012analysis}. Further, the components of \(\mathbf{\bar{u}}\) should be aranged in a SV-wise order:

\begin{equation}
    \begin{aligned}
        \mathbf{\bar{u}} =\left( \mathbf{\bar{u}_1},..., \mathbf{\bar{u}_N} \right) = \left(\bar{u}_{1,1}, \bar{u}_{1,2},...,\bar{u}_{1,k}, \bar{u}_{2,1},...,\bar{u}_{N,k}\right)
    \end{aligned}
\end{equation}
Subsequently, we proceed by calculating the reconstruction polynomials from the cell-averaged inital condition. The reconstruction is repetitive in every time step. Hence, we consider the reconstruction problem for the SV method in general: Given the cell-averaged state variables for all the control volumes in one spectral volume \(S_{i}\), construct a polynomial \(p_{i}\left(x\right)\) of degree at most \(k-1\) within \(S_{i}\) such that it is a smooth \(k\)th order accurate approximation to \(u\left(x \right)\) inside \( S_{i} \):

\begin{equation}
    \begin{aligned}
        p_{i} \left( x \right) = u \left( x \right) + \mathcal{O} \left(h^{k} \right), \quad i = 1,...,N.
    \end{aligned}
\end{equation}
Once we have the reconstruction polynomials, they can be employed to calculate approximations of \( u\left( x \right) \) at the CV boundaries.

\begin{equation}
    u_{i,j+1/2} := p_{i} \left( x_{i,j+1/2} \right) = u\left( x_{i,j+1/2} \right) + \mathcal{O} \left(h^{k} \right), \quad j = 0,...,k, \quad i = 1,...,N
\end{equation}
It is important to note, that we obtain two different values for \(u\) at the SV boundaries, because two different reconstruction polynomials intersect there.

The reconstruction problem has already been solved \cite[Problem 2.1]{cockburn1998essentially}. In problem 2.1 it is shown that for \( i =  1,...,N \), there exists a unique polynomial \( p_i \left(x \right)\) of degree at most \( k-1 \) in \(S_i \), whose cell averages in every control volume coincide with the given cell averages \( \bar{u}_{i,j}\)

\begin{equation}
    \bar{u}_{i,j} = \frac{1}{h_j} \int_{x_{i,j-1/2}}^{x_{i,j+1/2}} p_i \left( x \right) dx, \quad j = 1,...,k, \quad i = 1,...,N.
\end{equation}
Hence this reconstruction aligns with the conservation law, as the averaged integral of the reconstruction over each control volume reproduces the cell averages.

\begin{figure}[!h]
    \centering
    \includegraphics[width = \textwidth]{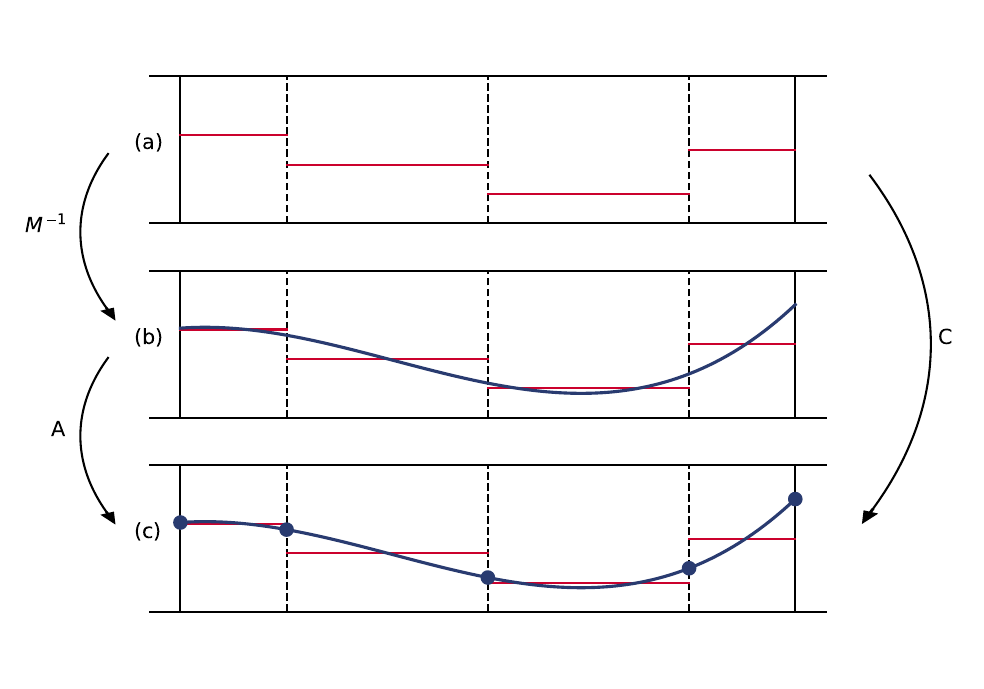}
    \caption{Polynomial reconstruction:
            (a) Cell averages, 
            (b) Reconstruction polynomial, 
            (c) Reconstructed boundary values.\\
             One can see the reconstruction within one spectral volume that is divided into four CVs. Initially cell averages are given for all four control volumes. In a first step the coefficients of the reconstruction polynomial are calculated via a multiplication of these cell averages by \(\mathbf{M^{-1}}\). Subsequently, the reconstruction polynomial is evaluated at all CV boundaries via a multiplication of the vector of coefficients by the matrix \( \mathbf{A}\). In practice, the two steps are of course combined in one multiplication of the given cell averages by the reconstruction matrix \( \mathbf{C}\).}
    \label{fig:reconstruction}
\end{figure}

We will now outline a simple way to construct \( p_i\left( x \right) \) numerically. The occuring matrices should have a low condition number. For this purpose, we will use the reliable Legendre polynomials \cite{isaacson2012analysis} for the reconstruction process. We express every polynomial \( p_i \left( x \right)\) as a combination of the first \( k\) Legendre polynomials \( L_m, \, m = 0,...,k-1\), using the inverse linear transformation of \eqref{eq:lintrans} to bring each spectral volume back into the interval \( \left[-1, 1 \right]\)

\begin{equation}
\begin{aligned}
     p_i \left(x \right) = l_{i,0} \cdot L_0 \left( \frac{2 \left( x-x_i \right)}{h} \right) + ... + l_{i,k-1} \cdot L_{k-1} \left( \frac{2\left(x - x_i \right)}{h} \right),\\
     \quad x \in S_i, \quad i = 1,...,N.
\end{aligned}   
\end{equation}
Hence,

\begin{equation}
    \begin{aligned}
       \frac{1}{h_j} \int_{x_{i,j-1/2}}^{x_{i,j+1/2}} p_i \left( x \right) dx &=
       \frac{1}{h_j} \int_{x_{i,j-1/2}}^{x_{i,j+1/2}} \sum_{m=0}^{k-1} l_{i,m} \cdot L_m \left( \frac{2 \left(x-x_i \right)}{h} \right)  dx\\[0.5cm]
       &=   \sum_{m=0}^{k-1} l_{i,m} \cdot \frac{1}{h_j} \int_{x_{i,j-1/2}}^{x_{i,j+1/2}} L_m \left( \frac{2 \left(x-x_i \right)}{h} \right)  dx.
    \end{aligned}
\end{equation}
The substitution \( s = \frac{2 \left(x - x_i \right)}{h}\) yields

\begin{equation}\label{eq:310}
    \begin{aligned}
       \bar{u}_{i,j} = \frac{1}{h_j} \int_{x_{i,j-1/2}}^{x_{i,j+1/2}} p_i \left( x \right) dx 
       =  \sum_{m=0}^{k-1} l_{i,m} \cdot \frac{1}{\tilde{h}_j} \int_{ \tilde{x}_{j-1/2}}^{ \tilde{x}_{j+1/2}} L_m \left( x \right) dx
    \end{aligned}
\end{equation}
where

\begin{align}
        \tilde{x}_{j+1/2} &= \frac{2 \left(x_{i,j+1/2} - x_i \right)}{h},\\[0.5cm]
        \tilde{h}_j &= \frac{h_j}{h} = \frac{2 \left(x_{i,j+1/2} - x_{i,j-1/2} \right)}{h} =  \tilde{x}_{j+1/2} - \tilde{x}_{j-1/2}
\end{align}
Examining \eqref{eq:310} for every \(p_i\), we obtain a system of linear equations to calculate the coefficients \(l_{i,m}, \, m = 0,...,k-1\). This system can be expressed using a \( k \times k\) matrix \( \mathbf{M} \), with entries defined as follows:

\begin{equation}
    M_{j,m} =  \frac{1}{\tilde{h}_j}  \int_{ \tilde{x}_{j-1/2}}^{ \tilde{x}_{j+1/2}} L_m \left( x \right) dx, \quad j = 1,...,k, \quad m = 0,...,k-1.
\end{equation}
The resulting linear system reads:

\begin{equation}
    \mathbf{\bar{u}_i} = \mathbf{M} \cdot \mathbf{l_i} \quad \Leftrightarrow \quad \mathbf{l_i} = \mathbf{M ^{-1}} \cdot \mathbf{\bar{u}_i},
\end{equation}
where for fixed \(i\) \( \mathbf{\bar{u}_i} \) is the vector of the \( (\bar{u}_{i,j})_{j=1}^{k}\) and \( \mathbf{l_i} \) is the vector of the coefficients \( (l_{i,m})_{m=0}^{k-1}\). To save computational cost, we calculate \( \mathbf{M ^{-1}} \) once and use it for every time step. It is worth noting that \( \mathbf{M}\) has to be non-singular, as the reconstruction polynomial \( p_i \left( x \right)\) is unique. Ultimately, we need to compute the approximate solution of \( u\left( x \right) \) at each CV boundary to update the state variable at the next time level:

\begin{equation}\label{eq:3.18}
\begin{aligned}
     u_{i,j+1/2} &= p_i \left( x_{i,j+1/2} \right) = \sum_{m = 0}^{k-1} l_{i,m} \cdot L_m \left(  \frac{2 \left(x_{i,j+1/2} - x_i \right)}{h} \right)\\
     &= \sum_{m = 0}^{k-1} l_{i,m} \cdot L_m \left( \tilde{x}_{j+1/2} \right).
\end{aligned}   
\end{equation}
Let \(\mathbf{u_i} \) be the vector of the reconstructed values \(u_{i,j+1/2}\). Equation \eqref{eq:3.18} can be written in a compact matrix vector form:

\begin{equation}
    \mathbf{u_i} = \mathbf{A} \cdot \mathbf{l_i}, \quad A_{j,m} = L_m \left( \tilde{x}_{j+1/2} \right), \quad j=0,...,k, \quad m = 0,...,k-1.
\end{equation}
Altogether, the approximate nodal value at the CV boundaries is given by

\begin{equation}\label{eq:316}
   \mathbf{u_i} = \mathbf{A} \cdot \mathbf{l_i} = \mathbf{A} \cdot \mathbf{M ^{-1}} \cdot \mathbf{\bar{u}_i}.
\end{equation}
Defining the reconstruction matrix as:

\begin{equation}
    \mathbf{C} = \mathbf{A} \cdot \mathbf{M ^{-1}}
\end{equation}
\eqref{eq:316} yields:
\begin{equation}
    \mathbf{u_i} = \mathbf{C} \cdot \mathbf{\bar{u}_i}.
\end{equation}
The reconstruction process is illustrated stepwise in figure \ref{fig:reconstruction}.

\clearpage

\subsection{The fundamental SV scheme}\label{sec:The fundamental SV scheme}
\begin{figure}[!h]
    \centering
    \includegraphics[width = \textwidth]{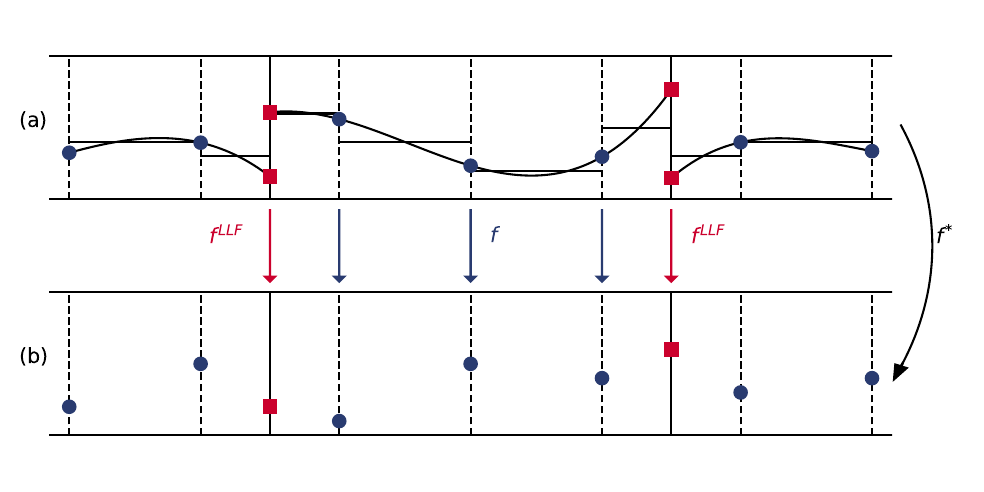}
    \caption{Numerical flux: (a) Reconstructed boudary values, (b) Numerical flux\\
    The computation of the numerical flux is depicted for one spectral volume that is divided into four CVs, with half of the adjacent SVs also visible. Initially, reconstructed values of \(u\) at the CV boundaries are provided, calculated from the cell averages \(\bar{u}_{i,j}\) as explained in the preceding chapter. At SV boundaries, two distinct values of \(u\) are obtained-one resulting from the reconstruction in the adjacent SV on the left and another from the reconstruction in the adjacent SV on the right. The ensuing Riemann problem is solved using an approximate Riemann solver. Consequently, at SV boundaries, the numerical flux is computed using the local Lax-Friedrichs flux \(f^{LLF} \) indicated by red squares.\\
    At inner CV boundaries, only one reconstructed value of \(u\) is obtained, as the reconstruction polynomial remains consistent throughout the entire spectral volume. Thus, at inner CV boundaries, the analytical flux \(f\) can be directly applied indicated by blue circles. In summary, this case-based definition distinctly determines the numerical flux \(f^{*}\) at every CV boundary.
    }
    \label{fig:flux}
\end{figure}

The fundamental spectral volume scheme presented in this section aims to provide a semi-discrete version of the given problem, allowing for numerical time integration. Analogue to the derivation of the FV method we start by integrating our differential equation \eqref{eq:conservation law} over one control volume \( C_{i,j} \):
\begin{equation}
    \begin{aligned}
         \int_{x_{i,j-1/2}}^{x_{i,j+1/2}} \partial_{t}u\left(x,t\right) +  \partial_{x} f \left(u\left(x,t\right)\right) dx = 0.
    \end{aligned}
\end{equation}
If \(u\) is sufficiently  smooth, the fundamental theorem of calculus yields:
\begin{equation}\label{eq:semidiscrete}
    \begin{aligned}
         h_j \cdot\frac{d \bar{u}_{i,j}}{dt} \left(t \right) + \left( f_{i,j+1/2} \left( t \right) - f_{i,j-1/2} \left( t \right)\right) = 0,
    \end{aligned}
\end{equation}
where
\begin{align}
        f_{i,j+1/2} \left( t \right) &= f \left(u\left(x_{i,j+1/2}, t\right) \right)\quad \text{and}\\[0.5cm]
        \frac{d \bar{u}_{i,j}}{dt} \left(t \right) &=  \frac{1}{h_j} \cdot \frac{d}{dt} \int_{x_{i,j-1/2}}^{x_{i,j+1/2}} u\left(x,t\right) dx.
\end{align}
Equation (\ref{eq:semidiscrete}) is equivalent to
\begin{equation}\label{eq:semidiscrete2}
    \frac{d \bar{u}_{i,j}}{dt} \left(t \right) = \frac{1}{h_j} \left( f_{i,j-1/2} \left( t \right) - f_{i,j+1/2} \left( t \right)\right).
\end{equation}
Like equation \eqref{eq:FV_semidiscrete} of the FV method, this equation can't be integrated directly, because we have no access to nodal values of \(u\) at SV boundaries. Therefore, we introduce a numerical flux \(f^*\) that approximates the flux at all cell volume boundaries
\begin{align}
    f_{i, j+1/2}^* \approx f(u(x_{i, j+1/2}, t)).
\end{align}
The concrete definition of \(f^*\) remains to be determined. As mentioned in the last section, we can use the analytical flux at the interior CV boundaries of a SV, due to the continuity of the reconstruction within every SV, i.e.
\begin{align}\label{eq:Anaflux}
    \begin{aligned}
         f_{i,j+1/2}^{*} = f\left(u_{i,j+1/2} \right), \quad j = 1,...,k-1.
    \end{aligned}
\end{align}
At the SV boundaries, on the other hand, we obtain two reconstructed values of the state variable. For example, at \(x_{i+1/2}\) the SV \(S_i\) adjoins \(S_{i+1}\). Hence one value \( u_{i,k+1/2} \) is due to the reconstruction in \(S_i\) and another value \( u_{i+1,1/2} \) due to the reconstruction in \(S_{i+1}\). In all numerical tests the local Lax-Friedrichs flux is used to calculate the numerical flux at SV boundaries. It is a simplification of the HLL flux introduced by Harten, Lax and van Leer \cite{harten1983upstream}.
\begin{equation}
    \begin{aligned}
&f^{HLL} \left(u_l, u_r \right) =
\begin{cases}
    f \left(u_l \right), & \quad a_l < 0 \\
    -\frac{a_l}{a_r - a_l} \cdot f \left( u_r \right) + \frac{a_r}{a_r - a_l} \cdot f \left( u_l \right) + \frac{a_l \cdot a_r}{a_r - a_l} \left(u_r - u_l \right), & \quad a_l < 0 < a_r\\
    f \left( u_r \right), & \quad a_r > 0
\end{cases}
    \end{aligned}
\end{equation}
where \( a_l\) and \(a_r\) are lower and upper bounds for the smallest and largest signal speeds respectively. The special choice of \(c_{max} = max\{ |a_l|, |a_r| \}\) and afterwards setting \(a_r = c_{max}, a_l = -c_{max}\) leads to the simpler local Lax-Friedrichs flux:

\begin{equation}\label{eq:LLF}
    \begin{aligned}
         f^{LLF}\left(u_l, u_r \right) = \frac{1}{2} \left( f\left( u_l \right) + f \left( u_r \right) \right) - c_{max} \cdot \frac{1}{2} \left(u_r - u_l \right).
    \end{aligned}
\end{equation}
To gain an intuitive understanding of the numerical flux, see figure \ref{fig:flux}.

For now, we assume periodic boundary conditions. Then by \eqref{eq:Anaflux} and \eqref{eq:LLF} the numerical flux is computed via

\begin{equation}\label{eq:fnum}
    \begin{aligned}
         &f_{i,j+1/2}^{*} =
\begin{cases}
    f^{LLF} \left( u_{N,k+1/2}, u_{1,1/2} \right) , &\quad \left(i = 1 \text{ and } j = 0 \right) \text{ or } \left( i = N \text{ and } j = k \right) \\
    f^{LLF} \left(u_{i,k+1/2}, u_{i+1,1/2} \right) , &\quad i = 1,...,N, \quad j = k \\
    f^{LLF} \left(u_{i-1,k+1/2}, u_{i,1/2} \right) , &\quad i = 1,...,N, \quad j = 0 \\
    f\left(u_{i,j+1/2}\right), &\quad i = 1,...,N, \quad j = 1,...,k-1 \\   
\end{cases}
    \end{aligned}
\end{equation}
where the first case is due to the periodic boundary condition. By the definition of the numerical flux \(f^{*}\) in \eqref{eq:fnum} the SV scheme is fully conservative.

Finally, we are able to state a scheme to compute the solution numerically. We use a third-order Runge-Kutta method \cite{shu1988total} for time integration of \eqref{eq:semidiscrete2}.
\begin{equation}\label{eq:Runge}
    \begin{aligned}
        \mathbf{\bar{u}^{\left(1\right)}} &= \mathbf{\bar{u}^{\left(0\right)}} + \Delta t \cdot \mathbf{D} \left( \mathbf{ \bar{u}^{\left(0\right)} }\right),\\
        \mathbf{\bar{u}^{\left(2\right)}} &= \frac{3}{4}  \mathbf{\bar{u}^{\left(0\right)}} + \frac{1}{4} \left( \mathbf{\bar{u}^{\left(1\right)}} + \Delta t \cdot \mathbf{D} \left( \mathbf{\bar{u}^{\left(1\right)}}\right) \right),\\
        \mathbf{\bar{u}^{\left(3\right)}} &= \frac{1}{3} \mathbf{\bar{u}^{\left(0\right)}} + \frac{2}{3}\left( \mathbf{\bar{u}^{\left(2\right)}} + \Delta t \cdot \mathbf{D}\left( \mathbf{\bar{u}^{\left(2\right)}} \right) \right),\\
    \end{aligned}
\end{equation}
where

\begin{equation}
         D_{i,j}\left( \mathbf{\bar{u}} \right) = \frac{1}{h_j} \left(f_{i,j-1/2}^{*} - f_{i,j+1/2}^{*} \right).
\end{equation}
is computed by formula \eqref{eq:fnum}.

\subsection{First numerical Test}\label{sec:First numerical test}
\begin{figure}[!h]
    \centering
    \includegraphics[width=0.65\linewidth]{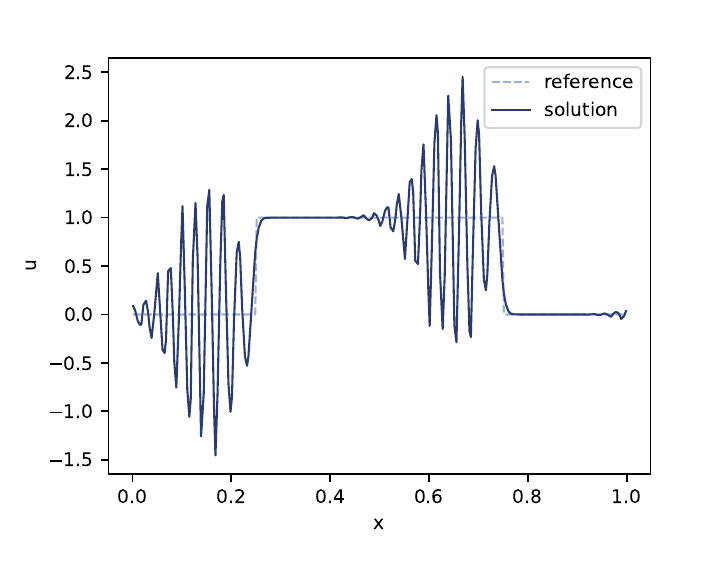}
    \caption{Solution at \(t = 1.0\)}
    \label{fig:rectangle_unstable}
\end{figure}

Now, we perform some simple test to analyse the basic SV implementation. To ensure comparability, we adopt a problem proposed by Wang \cite{wang2002spectral}. Consider the linear advection equation defined by the flux function \(f\left( u \right) = v \cdot u\). Thus \eqref{eq:conservation law} reduces to
\begin{align}
    \partial_t u \left(x, t \right) + v \cdot \partial_x u\left(x, t \right) = 0, \quad u\left(x, 0\right) = u_0\left(x \right).
\end{align}
We choose the spatial domain \( \left[0,1\right]\) and the velocity \(v = 1\). Periodic boundary conditions are assumed and the initial condition is given by
\begin{equation}
    u_{0}\left(x\right) = 
    \begin{cases}
        1, \quad \frac{1}{4} \leq x \leq \frac{3}{4}\\
        0, \quad \text{otherwise.}
    \end{cases}
\end{equation}
We use \(N = 60\) spectral volumes and \(k = 4\) control volumes each for the numerical computation. The numerical solution is depicted in figure \ref{fig:rectangle_unstable} alongside an analytical reference at \( t = 1.0\). Which means that the rectangle has passed one complete round through the domain.

\begin{figure}[!h]
    \centering
    \includegraphics[width=0.65\linewidth]{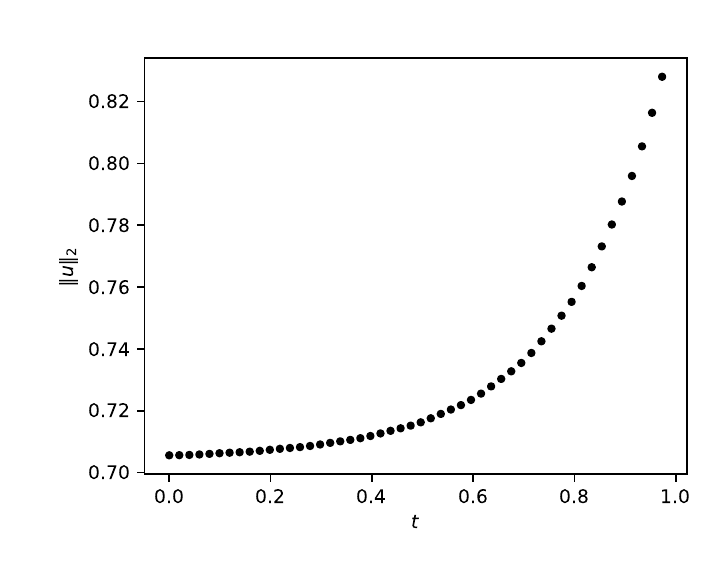}
    \caption{\(L^{2}\)-norm of the numerical solution over time.}
    \label{fig:rectangle_unstable_L2}
\end{figure}
\clearpage
Evidently, the solution developed pronounced oscillations, even though the rectangle has passed the domain only once. This behavior indicates a clear instability in the scheme. To verify this assumption, the \(L^{2}\)-norm of the numerical solution was calculated for several points in time, and the results are shown in figure \ref{fig:rectangle_unstable_L2}.
The distinct slope observed in the \(L^{2}\)-norm over time validates the suspected instability \cite{svard2014review}. Consequently, the scheme presented in the previous section cannot always be directly applied for practical purposes, and a method to stabilize our Spectral-Volume scheme is needed.

\section{Entropy Stabilization}\label{sec:Entropy Stabilization}

In this chapter, we provide an overview of an entropy rate stabilization approach for the Discontinuous Galerkin (DG) method as presented by Klein \cite{klein2023stabilizing2} and demonstrate its adaption to the Spectral-Volume scheme.
After grasping the basic idea of the stabilization method in section \ref{sec:Basic idea of the Stabilization method}, we start by an estimate for the highest possible entropy dissipation speed in section \ref{sec:Estimate for the highest possible entropy dissipation}. In the following chapter a direction for the dissipation is constructed in a physically motivated manner based on the heat equation. The correction size for the required entropy dissipation will be calculated in section \ref{sec:Correction size}. Finally, in section \ref{sec:Stabilized SV scheme}, we present the stabilized SV scheme. During this chapter, we exclusively use cell averages. Hence, we no longer employ the upper bar introduced in the previous chapter. Still all vectors of cell averages are presented in bold typeface, as are all matrices.

\subsection{Basic Idea}\label{sec:Basic idea of the Stabilization method}

As we use cell averages for our spectral volume method, we will in general approximate an integral over \(S_{i}\) by
\begin{equation}\label{eq:CV quadrature}
    \int_{S_{i}} g\left(u\right) \, dx \approx \sum_{j=1}^{k} h_{j} g\left(u_{i,j} \right),
\end{equation}
where \(g: \mathbb{R}^{m} \rightarrow \mathbb{R}\) is a given function, for example, an entropy \(U\). Similarly, a discretisation of the inner product will be given by

\begin{equation}\label{eq:CV scalarproduct}
     \int_{S_{i}} f\left(u\right)g\left(u\right) \, dx \approx \sum_{j=1}^{k} h_{j} g\left(u_{i,j} \right) f\left(u_{i,j} \right)
     =: \left \langle f, g \right\rangle_{S_{i}}.
\end{equation}
Thus, for a numerical approximation of the total entropy in one SV \(S_{i}\) at time \(t\), we obtain

\begin{equation}
    E_{u,i}\left(t\right) = \int_{x_{i-1/2}}^{x_{i+1/2}} U\left(u\left(x, t\right)\right) \, dx 
    \approx \sum_{j=1}^{k} h_{j} U\left(u_{i,j}\left(t\right)\right).
\end{equation}
Note, that the discretizations for the integral and scalar product on one SV are not exact. This can be seen as follows:\\
The approximations \eqref{eq:CV quadrature} and \eqref{eq:CV scalarproduct} are obviously exact if \(u = const\). On the other hand if we consider \(u\) to be linear, for example \(u = x\) and \(g\left(u\right) = u^{2}\), the approximation is not exact as can be shown by a simple calculation:

\begin{align}
    \int_{S_{i}} g\left(u\right) \, dx = \int_{x_{i-1/2}}^{x_{i+1/2}} x^{2} \, dx = \frac{1}{3} \left(x_{i+1/2}^{3} - x_{i-1/2}^{3} \right)
\end{align}
but
\begin{align}
    \sum_{j=1}^{k} h_{j} g\left(u_{i,j} \right) = \sum_{j=1}^{k} \left( x_{i,j+1/2} - x_{i,j-1/2} \right) \left( \frac{x_{i,j-1/2} + x_{i,j+1/2}}{2} \right)^{2}.
\end{align}
For convenience, we assume \(x_{i-1/2} = -1, \, x_{i+1/2} = 1\) and four equidistant control volumes. Then

\begin{align}
    \int_{S_{i}} g\left(u\right) \, dx  = \frac{1}{3} \left(x_{i+1/2}^{3} - x_{i-1/2}^{3} \right) = \frac{2}{3} \neq 0.625 = \sum_{j=1}^{k} h_{j} g\left(u_{i,j} \right). 
\end{align}
Now, we have all necessary definitions to start our exploration of the stabilization method. The procedure is as follows:

\begin{enumerate}
    \item Calculate the time derivatives \(\frac{d u_{i,j}}{dt}\) of all cell averages.
    \item Estimate the highest possible entropy dissipation speed \(\sigma_{i + 1/2}\) at the interface of \(S_{i}\) and \(S_{i+1}\).
    \item Compute a correction direction \(v_{i}\).
    \item Determine the necessary size \(\lambda_{i}\) of the correction to achieve that

    \begin{equation}
        \frac{d \tilde{u}_{i,j}}{dt} = \frac{d u_{i,j}}{dt} + \lambda_{i} v_{i}
    \end{equation}

    satisfies

    \begin{equation}
        \frac{d E_{\tilde{u},i}}{dt} = \left \langle \frac{dU}{du}, \frac{d u}{dt} + \lambda_{i} v_{i} \right\rangle_{S_{i}} \leq \sigma_{i} + F_{i-1/2}^{*} - F_{i+1/2}^{*},
    \end{equation}

    where \(F^{*}\) is a numerical entropy fux. In hope for an entropy inequality, we employ the (local) Lax-Friedrichs flux for \(F^{*}\) as well \cite{tadmor1984large, tadmor1984numerical}, i.e.

    \begin{equation}
        F^{*}_{i+1/2} = \frac{1}{2} \left( F\left( u_{i,k} \right) + F \left( u_{i+1,1} \right) \right) - c_{max} \cdot \frac{1}{2} \left(u_{i+1,1} - u_{i,k} \right).
    \end{equation}
\end{enumerate}
For the first step, equation \eqref{eq:semidiscrete2} from the previous chapter is applicable. In chapter \ref{sec:Basic idea of the Stabilization method} it was established that entropy dissipation can only occur in the absence of smooth solutions, as the entropy equality \eqref{eq:entropy equality} holds otherwise. Consequently, dissipation does not transpire within the spectral volumes, given the smoothness of the reconstruction polynomials. Instead, dissipation is dictated by the SV boundaries where distinct reconstruction polynomials intersect.

\subsection{Estimate for the highest possible entropy dissipation}\label{sec:Estimate for the highest possible entropy dissipation}

Given the occurence of local Riemann problems at SV boundaries in every time step, we investigate a general Riemann problem with an initial condition given by:

\begin{equation}
    u_{0}\left( x \right) =
    \begin{cases}
        u_{l}, \quad x < 0,\\
        u_{r}, \quad x > 0.
    \end{cases}
\end{equation}
From \cite[Lemma 2]{klein2023stabilizing2} we obtain the following theorem.\\

\begin{theorem}
    Given bounds \(a_{l}, a_{r}\) on the fastest signal speed to the left and to the right, respectively. Let \(M \geq \max \{ | a_{l} |, |a_{r}| \} \). The entropy dissipation speed of a Riemann problem solution on the interval \(\theta = \left] -M, M\right[\) is bounded from below by

    \begin{equation}\label{eq:entropy dissipation rate estimate}
        \sigma^{\theta} \geq \left(a_{r} - a_{l} \right) U\left(u_{lr} \right) + a_{l} U\left(u_{l}\right) - a_{r} U\left(u_{r} \right) + F\left(u_{l}\right) - F \left( u_{r} \right),
    \end{equation}
    where

    \begin{equation}
        u_{lr} = \frac{a_{r} u_{r} - a_{l} u_{l} + f\left(u_{l} \right) - f\left(u_{r} \right)}{a_{r} - a_{l}}.
    \end{equation}
    The proof can be found in the literature \cite[Lemma 2]{klein2023stabilizing2}.
\end{theorem}
To align with the used local Lax-Friedrichs flux, we should adapt this inequality by choosing \(c_{max} = \max\{ |a_l|, |a_r| \}\) and subsequaently updating \(a_r = c_{max}, a_l = -c_{max}\). This leads to

\begin{equation}
    u_{lr} =  \frac{u_{l} + u_{r}}{2} + \frac{f\left(u_{l} \right) - f\left(u_{r} \right)}{2c_{max}}
\end{equation}
and \eqref{eq:entropy dissipation rate estimate} simplifies to
\begin{equation}\label{eq:entropy estimate}
    \sigma^{\theta} \geq c_{max}\left(2 U\left(u_{lr} \right) - U\left(u_{l}\right) - U\left(u_{r} \right) \right) + F\left(u_{l}\right) - F \left( u_{r} \right).
\end{equation}
This lower bound will inevitably overstate the highest possible entropy dissipation speed, given its negative nature. However, aiming for the most dissipative weak solution, we will employ this lower bound as an approximation of the highest possible entropy dissipation speed \cite{tadmor1984large, tadmor1984numerical}.\\
To incorporate the entropy dissipation speed approximation in our SV scheme, we apply \eqref{eq:entropy estimate} on all local Riemann problems given by the reconstructed values of \(u\) at the SV boundaries. This leads to
\begin{equation}
\begin{aligned}\label{eq:sigma}
     \sigma_{i+1/2} = \, &c_{max}\left(2 U\left(u_{m} \right) - U\left(u_{i,k+1/2}\right) - U\left(u_{i+1,1/2} \right) \right)\\
    &+ F\left(u_{i,k+1/2}\right) - F \left( u_{i+1,1/2} \right)
\end{aligned}   
\end{equation}
with
\begin{equation}
    u_{m} = \frac{u_{i,k+1/2} + u_{i+1,1/2}}{2} + \frac{f\left(u_{i,k+1/2} \right) - f\left(u_{i+1,1/2} \right)}{2c_{max}}
\end{equation}
using periodic adjustments at the boundaries if necessary.

\subsection{Entropy dissipation direction}\label{sec:Entropy dissipation direction}

Our objective now is to determine a direction to correct the time derivative of the SV scheme, ensuring entropy dissipation.
In accordance with the DG version \cite{klein2023stabilizing2} the correction direction will be based on the construction of operators that can regularize a solution \(u\). These operators, termed filters, take the form of special Hilbert-Schmidt operators \(K\).\\

\begin{definition}\label{def:filter}
    \textnormal{(Filter):} An operator \(K: L^{2} \left( \Omega \right) \rightarrow L^{2} \left( \Omega \right) \) is designated as filter if it is an integral operator satisfying

    \begin{equation}
        \left[Ku\right]\left(x\right) = \int_{\Omega} k\left(x, y\right) u\left(y\right) \, dy, \quad \text{with} \quad \forall x \in \Omega: 
        \int_{\Omega} k\left(x,y \right) \, dy = 1.
    \end{equation}
    The kernel \(k\) is required to have bounded Hilbert-Schmidt norm. The pointwise evaluation of \(\left[Ku\right]\) conforms to a weighted average of \(u\).
\end{definition}
Of particular interest are conservative filters, ensuring the preservation of our basic SV scheme's conservation property when applied to cell averages. Fubini's theorem provides:

\begin{theorem}\label{th:conservative filter}
    \textnormal{(Conservative filter):}
    A filter \(K: L^{2} \left( \Omega \right) \rightarrow L^{2} \left( \Omega \right) \) is conservative, i.e.

    \begin{equation}
        \int_{\Omega} \left[Ku\right] \left(x \right) \, dx = \int_{\Omega} u\left(x\right) \, dx,
    \end{equation}
    if the corresponding kernel \(k\) has mass one, i.e.

    \begin{equation}
        \forall y \in \Omega: \quad \int_{\Omega} k\left(x,y \right) \, dx = 1.
    \end{equation}
    For the complete proof, refer to \cite[Lemma 4]{klein2023stabilizing2}.
\end{theorem}
Our aspiration is for the filter \(K\) to impart a smoothing effect and to dissipate entropy. The selection for an appropriate filter \(K\) is guided by the following theorem:

\begin{theorem}\label{th:dissipative filter}
    \textnormal{(Dissipative filter):}
    A filter \(K\) is deemed dissipative if for all convex entropy functions \(U\)
    \begin{equation}
        E_{Ku} = \int_{\Omega} U\left( \left[Ku\right]\left(x\right)\right) \, dx \leq \int_{\Omega} U\left(u\left(x\right) \right) \, dx = E_{u}
    \end{equation}
    holds. For a conservative filter the above inequality holds if the corresponding kernel is non-negative, i.e.
    \begin{align}\label{eq:positive kernel}
        \forall x,y \in \Omega: k\left(x,y\right) \geq 0.
    \end{align}
    If the corresponding kernel \(k\) of a filter \(K\) satisfies (\ref{eq:positive kernel}) we refer to the filter as positive. The proof can be found in the literature \cite[Theorem 1]{klein2023stabilizing2}.
\end{theorem}
Evidently, the continuous filter cannot be applied directly on our numerical SV scheme. The next step involves modifying these concepts for discrete application. The discrete version of the continuous conservative, positive and hence dissipative filter is defined as follows:

\begin{definition}\label{def:discrete filter}
    \textnormal{(Discrete conservative and positive filter):}
    We designate a matrix \(\mathbf{Y} \in \mathbb{R}^{k \times k}\) as a filter if
    \begin{equation}
        \forall j \in \{1,...,k\}: \quad \sum_{l=1}^{k} Y_{jl} = 1
    \end{equation}
    holds by analogy with definition \ref{def:filter}. It is said to be conservative, if
    \begin{equation}\label{eq:discrete conservation}
        \forall l \in \{1,...,k\}: \quad \sum_{j=1}^{k} h_{j} Y_{jl} = h_{l}
    \end{equation}
    is satisfied, see theorem \ref{th:conservative filter}. Further we call it positive if
    \begin{equation}
        \forall j,l \in \{1,...,k\}: \quad Y_{jl} \geq 0
    \end{equation}
    holds corresponding to theorem \ref{th:dissipative filter}.
\end{definition}
This definition may seem at odds with the continuous version of the filter \(K\), especially if we consider the matrix \(\mathbf{Y}\) to be the discrete representation of the integral kernel \(k\). The contradiction is resolved by considering the relation

\begin{equation}
    Y_{jl} = h_{l} \tilde{Y}_{jl}
\end{equation}
and interpreting \(\mathbf{\tilde{Y}}\) as the discrete version of the integral kernel \(k\). Consequently, the integral operator \(K\) simplifies to a mere matrix-vector multiplication in the discrete case:
\begin{equation}
    \left[Ku\right]\left(x\right) = \int_{\Omega} k\left(x, y\right) u\left(y\right) \, dy \quad \leftrightarrow \quad
    v_{i,j} = \sum_{l=1}^{k} h_{l} \tilde{Y}_{jl} u_{i,l} = \sum_{l=1}^{k} Y_{jl} u_{i,l}.
\end{equation}
The relation (\ref{eq:discrete conservation}) arises from
\begin{align*}
    \int_{S_{i}} \left[Ku \right] \left(x \right) \, dx
     \approx \sum_{j=1}^{k} h_{j} \sum_{l=1}^{k} Y_{jl} u_{i,l}
     = \sum_{l=1}^{k} u_{i,l} \sum_{j=1}^{k} h_{j} Y_{jl}
     = \sum_{l=1}^{k} h_{l} u_{i,l}
     = \int_{S_{i}} u \left( x \right) \, dx
\end{align*}

\begin{theorem}
    Analogous to theorem \ref{th:conservative filter}, we establish that a discrete filter \(Y\) that is positive and conservative is also dissipative in the sense that
    \begin{equation}
         \sum_{j=1}^{k} h_j U\left( v_{i,j} \right) \leq \sum_{l=1}^{k} h_l U\left( u_{i,l} \right),
         \quad v_{i,j} = \sum_{l=1}^{k} Y_{j,l} u_{i,l}.
    \end{equation}
    That means, a filter is dissipative if Jensen's inequality holds in a discrete sense.
    
    Proof: By Jensen's inequality and definition \ref{def:discrete filter} we get
    \begin{align*}
        \sum_{j=1}^{k} h_j U\left( v_{i,j} \right) 
        = \sum_{j=1}^{k} h_j U\left( \sum_{l=1}^{k} Y_{j,l} u_{i,l} \right) 
        \leq \sum_{j=1}^{k} \sum_{l=1}^{k}  h_j Y_{j,l} U\left( u_{i,l} \right).
    \end{align*}
    The result follows by using again definition \ref{def:discrete filter}
    \begin{align*}
        \sum_{j=1}^{k} \sum_{l=1}^{k}  h_j Y_{j,l} U\left( u_{i,l} \right)
        = \sum_{l=1}^{k} U\left( u_{i,l} \right) \left( \sum_{j=1}^{k}  h_j Y_{j,l} \right)
        = \sum_{l=1}^{k} h_{l} U\left( u_{i,l} \right).
    \end{align*} 
\end{theorem}
As we will see later, we can construct positive, conservative filters using so-called filter generators:

\begin{definition}
    \textnormal{(Conservative and positive filter generator):}
    Let \(G \in \mathbb{R}^{k \times k}\) a square matrix. we call \(G\) a filter generator if
    \begin{equation}
        \forall j \in \{1,...,k\}: \sum_{l=1}^{k} G_{jl} = 0
    \end{equation}
    holds. It will be conservative if
    \begin{equation}
         \forall l \in \{1,...,k\}: \sum_{j=1}^{k} h_{j} G_{jl} = 0
    \end{equation}
    is satisfied. Further we call it positive, if
    \begin{equation}
        \forall l \in \{1,...,k\}, \quad \forall j \in \{1,...,l-1,l+1,...,k\}: \quad G_{jl} \geq 0
    \end{equation}
    holds.
\end{definition}
Generators and filters have the following connection \cite[Lemma 5]{klein2023stabilizing2}:

\begin{theorem}\label{th:generators and filters}
    \textnormal{(Generators and filters):}
    \begin{equation}
        \mathbf{G} \text{ conservative as generator} \quad \Rightarrow \quad \mathbf{Y} = \mathbf{I} + \Delta t \mathbf{G} \text{ conservative as filter.}
    \end{equation}
    Let further \(\Delta t \max_{l} |G_{ll}| \leq 1\). Then it follows

    \begin{equation}
         \mathbf{G} \text{ positive as generator}  \quad \Rightarrow \quad \mathbf{Y} \text{ positive as filter.}
    \end{equation}
    The proof can be found in the literature \cite[Lemma 5]{klein2023stabilizing2}.
    
\end{theorem}
As in \cite{klein2023stabilizing2}, we will construct a filter generator based on the heat equation, thus possessing a direct physical interpretation. However, applying the procedure for the DG method directly to the SV scheme is not feasible as the DG procedure relies on the nodal basis of the DG method, whereas the SV approach is based on cell averages. Therefore, as an initial step, we introduce a finite volume discretization of the heat equation.

\vspace{0.5cm}

\textbf{Finite volume discretization of the heat equation:}\\
Consider the heat equation
\begin{equation}\label{eq:heat equation}
    \frac{\partial u}{\partial t} = \Delta u
\end{equation}
for an arbitrary spectral volume \(S_{i}\) in conjunction with homogenous Neumann boundary conditions. As the index \(i\) does not matter for our derivation, we will denote \(u_{j} = u_{i,j}\) for brevity. The heat equation \eqref{eq:heat equation} can be rewritten as conservation law
\begin{equation}
    \frac{\partial u}{\partial t} + \frac{\partial f}{\partial x} = 0
\end{equation}
with flux function

\begin{equation}
    f\left( \nabla u\right) = - \nabla u.
\end{equation}
Therefore, we define the discrete heat flux \(f^{h}\) as

\begin{equation}
    f^{h}_{0+1/2} = f^{h}_{k+1/2} = 0, \quad f^{h}_{j+1/2} = -2 \cdot \frac{u_{j+1} - u_{j}}{h_{j+1} + h_{j}}, \quad j = 1,...,k-1.
\end{equation}
An semidiscrete version of the heat equation is then given by

\begin{equation}\label{eq:heat equation semidiscrete}
    \frac{d u_j}{dt} = \frac{f^{h}_{j-1/2} - f^{h}_{j+1/2}}{h_{j}}, \quad j = 1,...,k.
\end{equation}
This can be expressed in matrix-vector form:

\begin{equation}\label{eq:H}
    \frac{d\mathbf{u}}{dt} = \mathbf{H} \mathbf{u}.
\end{equation}

\begin{theorem}
    The matrix \( \mathbf{H}\) defined by \eqref{eq:heat equation semidiscrete} and \eqref{eq:H} satisfies all requirements of a conservative, positive filter generator. Therefore, it will be used for all numerical computations and we denote \(\mathbf{H} = \mathbf{G}\).
    
    Proof: The matrix \( \mathbf{H}\) is by construction a tridiagonal matrix and the non-zeros entries are given by:
    \begin{equation}\label{eq:H boundary entries}
        \begin{aligned}
            H_{11} &= - \frac{2}{h_{1} \left( h_{1} + h_{2} \right)}, \quad & &H_{12} = \frac{2}{h_{1} \left( h_{1} + h_{2} \right)},\\[0.2cm]
            H_{kk-1} &= \frac{2}{h_{k} \left( h_{k-1} + h_{k} \right)}, \quad & &H_{kk} = - \frac{2}{h_{k} \left( h_{k-1} + h_{k} \right)},
        \end{aligned}
    \end{equation}
    and for \(j = 2,...,k-1\) we have
    \begin{equation}\label{eq:H entries}
        \begin{aligned}
            H_{jj-1} &= \frac{2}{h_{j} \left( h_{j-1} + h_{j} \right)},\quad H_{jj+1} = \frac{2}{h_{j} \left( h_{j} + h_{j+1} \right)},\\[0.2cm]
        H_{jj} &= - \frac{2}{h_{j}} \left( \frac{1}{h_{j-1} + h_{j}} + \frac{1}{h_{j} + h_{j+1}} \right).
        \end{aligned}
    \end{equation}
    Thus, the matrix \( \mathbf{H}\) satisfies
    \begin{equation}
        \forall j \in \{1,...,k\}: \sum_{l=1}^{k} H_{jl} = 0
    \end{equation}
    which means \( \mathbf{H}\) is a filter generator. Secondly, \(\mathbf{H}\) is conservative, because \( \forall l \in \{2,...,k-1\}:\)
    \begin{align*}
         \sum_{j=1}^{k} h_{j} H_{jl} &= h_{l-1}H_{l-1l} + h_{l}H_{ll} + h_{l+1}H_{l+1l}\\
        &= h_{l-1} \cdot \frac{2}{h_{l-1} \left( h_{l-1} + h_{l} \right)} - h_{l} \cdot \frac{2}{h_l} \left( \frac{1}{h_{l-1} + h_{l}} + \frac{1}{h_{l} + h_{l+1}} \right)\\
        &\quad+ h_{l+1} \cdot \frac{2}{h_{l+1} \left(h_{l} + h_{l+1} \right)}\\[0.2cm]
        &= 2 \left(\frac{1}{ h_{l-1} + h_{l} } -  \left(\frac{1}{h_{l-1} + h_{l}} + \frac{1}{h_{l} + h_{l+1}} \right) +  \frac{1}{h_{l} + h_{l+1}} \right)\\
        &= 0
    \end{align*}
    and the same is true for \(l = 1\) and \(l = k\). Finally by \eqref{eq:H boundary entries} and \eqref{eq:H entries} \(\mathbf{H}\) is positive as
    \begin{equation}
        \forall l \in \{1,...,k\}, \quad \forall j \in \{1,...,l-1,l+1,...,k\}: \quad H_{jl} \geq 0
    \end{equation}
    holds.
\end{theorem}
Note, that we have no problems with the positivity of \(\mathbf{H}\), although this was an issue with the construction of a suitable \(\mathbf{G}\) for the DG version. The next theorem shows, that the matrix \(\mathbf{H}\) defined above dissipates everything that is not constant.

\begin{theorem}\label{th:nullspace of H}
    The nullspace of the matrix \(\mathbf{H}\) defined by \eqref{eq:H} consists only of constants.
    
    Proof: By the filter property, we know that all constants are contained in the nullspace of \(\mathbf{H}\). To check whether the constant vectors are the only elements of the nullspace of \( \mathbf{H} \), we use the Gauss algorithm. For brevity let \(r_{i}\) be the \(i\)th row of \( \mathbf{H}\).
    
    \textbf{Gauss algorithm:}
    \begin{enumerate}
        \item \(r_{2} = r_{2} + \frac{h_{1}}{h_{2}} \cdot r_{1}\)
        \item \(r_{3} = r_{3} + \frac{h_{2}}{h_{3}} \cdot r_{2}\)
        \item ...
        \item \(r_{k} = r_{k} + \frac{h_{k-1}}{h_{k}} \cdot r_{k-1}\)
        \item Multiply all rows by \(-1\).
    \end{enumerate}
    The resulting matrix has the form:
    \begin{equation}
        \mathbf{\tilde{A}} = 
        \begin{bmatrix}
            1 & -1 & 0 & 0 & \cdots & 0 \\
            0 & 1 & -1 & 0 & \cdots & 0 \\
            \vdots & \ddots & \ddots & \ddots & \ddots & \vdots \\
            0 & \cdots & 0 & 1 & -1 & 0 \\
            0 & \cdots & 0 & 0 & 1 & -1 \\
            0 & \cdots & \cdots & 0 & 0 & 0 \\
        \end{bmatrix}
    \end{equation}
    Next we can bring the matrix into reduced row-echelon form by:
    
    \begin{enumerate}
        \item \(r_{k-2} = r_{k-2} + r_{k-1}\)
        \item \(r_{k-3} = r_{k-3} + r_{k-2}\)
        \item ...
        \item \(r_{1} = r_{1} + r_{2}\)
    \end{enumerate}
    The reduced row-echelon form of \(\mathbf{H} \) is then given by:

    \begin{equation}
        \mathbf{A} = 
        \begin{bmatrix}
            1 & 0 & 0 & \cdots & 0 & -1 \\
            0 & 1 & 0 & \cdots & 0 & -1 \\
            \vdots & \ddots & \ddots & \ddots & \vdots & \vdots \\
            0 & \cdots & 0 & 1 & 0 & -1 \\
            0 & \cdots & 0 & 0 & 1 & -1 \\
            0 & \cdots & \cdots & 0 & 0 & 0 \\
        \end{bmatrix}
    \end{equation}
    Therefore the nullspace of \( \mathbf{H}\) actually consists only of constants.
\end{theorem}
Using theorem \ref{th:nullspace of H}, we can show that the chosen direction indeed dissipates entropy \cite[Lemma 6]{klein2023stabilizing2}:

\begin{theorem}
    If \(\mathbf{u}\) is non-constant with admissable strictly convex entropy \(U\) and the nullspace of \(\mathbf{H}\) consists only of constants, then
    \begin{equation}
        \left \langle \frac{dU}{du}, \mathbf{H} \mathbf{u_{i}} \right \rangle_{S_{i}} < 0.
    \end{equation}
    If \(U\) is merely convex, only \( \leq \) applies in the inequality above.
    Keeping theorem \ref{th:nullspace of H} in mind the proof can be found in the literature \cite[Lemma 6]{klein2023stabilizing2}.
\end{theorem}
Note, that we don't have to assume that the nullspace of \( \mathbf{H} = \mathbf{G}\) consists only of constants as in the DG case \cite[Lemma 6]{klein2023stabilizing2}, because of theorem \ref{th:nullspace of H}.

\subsection{Correction size}\label{sec:Correction size}
After calculating the suitable direction \(\mathbf{v_{i}} = \mathbf{Hu_{i}} = \mathbf{Gu_{i}}\), we continue with
the procedure outlined for the the DG version \cite[chapter 3.1]{klein2023stabilizing2} to compute an appropriate value of the correction size \(\lambda_i\) in two steps. First, we calculate

\begin{equation}\label{eq:lambda_ED}
    \lambda_{ED}^{i} = \max \left( 0, -\frac{\left \langle \frac{dU}{du}, \frac{du}{dt} \right \rangle_{S_{i}} - \left(F^*_{i-1/2} - F^{*}_{i+1/2} \right)}{\left \langle \frac{dU}{du}, \mathbf{v_{i}} \right \rangle_{S_{i}}} \right)
\end{equation}
to enforce the per spectral volume entropy dissipativity:

\begin{equation}
    \left \langle  \frac{dU}{du}, \frac{du}{dt} + \lambda^{i}_{ED} \mathbf{v_{i}} \right \rangle_{S_{i}} \leq F^{*}_{i-1/2} - F^{*}_{i+1/2},
\end{equation}
In the second step, we enforce a suffiently high enough entropy rate. This is achieved by splitting the estimated entropy dissipation rate of the adjacent Riemann problems. In general, every Riemann problem is connected with an approximation of the highest possible entropy rate \(\sigma^{\theta}\). The entropy rate approximations are then split between the two adjecent spectral volumes, leading to two additional parts of the correction size, one from the left and another from the right adjacent Riemann problem of the SV under consideration.

\begin{equation}\label{eq:lambda_ER}
    \begin{aligned}
        \lambda^{i}_{ER,l} = \max \left( 0, \frac{\sigma_{i-1/2}}{ \left \langle \frac{dU}{du}, \mathbf{v_{i-1}} \right \rangle_{S_{i-1}} + 
        \left \langle \frac{dU}{du}, \mathbf{v_{i}} \right \rangle_{S_{i}}  } \right),\\
        \lambda^{i}_{ER,r} = \max \left( 0, \frac{\sigma_{i+1/2}}{ \left \langle \frac{dU}{du}, \mathbf{v_{i}} \right \rangle_{S_{i}} + 
        \left \langle \frac{dU}{du}, \mathbf{v_{i+1}} \right \rangle_{S_{i+1}}  } \right).
    \end{aligned}
\end{equation}
All parts are then added together:

\begin{equation}\label{eq:lambda_sum}
    \lambda^{i}_{sum} = \lambda^{i}_{ED} + \lambda^{i}_{ER,l} + \lambda^{i}_{ER,r}.
\end{equation}
In a last step, we introduce an upper limit \( \lambda_{max}\) and compute 

\begin{equation}\label{eq:lambda_final}
    \lambda^{i} = \min \left( \lambda_{max}, \lambda^{i}_{sum} \right)
\end{equation}
to enforce the limit dictated by theorem \ref{th:generators and filters}. For a clarification of the notation used, note that
\begin{align}
    \left \langle \frac{dU}{du}, \frac{du}{dt} \right \rangle_{S_{i}} &= \sum_{j=1}^{k}h_{j} \frac{dU}{du} \left(u_{i,j} \right) \cdot \frac{f_{i, j-1/2}^{*} - f_{i, j+1/2}^{*}}{h_{j}},\\[0.5cm]
    \left \langle \frac{dU}{du}, \mathbf{v_{i}} \right \rangle_{S_{i}} &= \sum_{j=1}^{k}h_{j} \frac{dU}{du} \left(u_{i,j} \right) \cdot v_{i,j},
\end{align}
where
\begin{equation}
    \mathbf{v_{i}} = \mathbf{G u_{i}}, \quad \mathbf{v_{i}} = \left( v_{i,1},...,v_{i,k} \right).
\end{equation}

\subsection{The stabilized SV scheme}\label{sec:Stabilized SV scheme}

Now, the complete SV scheme is stated. Steps 1 to 5 of the basic SV scheme remain unchanged, but a step to calculate the needed filter generator \(\mathbf{G}\) from the heat equation is added. We still use the third order Runge-Kutta method in step 6. Only the routine for the Euler step has to be modified.

\textbf{Adapted euler step:} (input: \(\mathbf{\bar{u}^{n}}\), output: \(\mathbf{\bar{u}^{n+1}}\) )
\begin{enumerate}
    \item From the cell averages at time t, calculate the reconstructed values of u at the CV boundaries using the reconstruction matrix \(\mathbf{C}\):
    \begin{equation}
        \mathbf{u_i^n} = \mathbf{C} \cdot \mathbf{\bar{u}_i^n}.
    \end{equation}

    \item Compute the numerical flux \( f_{i,j+1/2}^{*,n} \) at the CV boundaries using formula \eqref{eq:fnum}.

    \item For all local Riemann-problems occuring at SV boundaries calculate the entropy dissipation speed approximation \( \sigma_{i+1/2}^{n} \) using \eqref{eq:sigma}.

    \item Compute the LLF numerical entropy flux \(F^{*,n}_{i+1/2}\) at all SV boundaries.

    \item Use formulas \eqref{eq:lambda_ED}, \eqref{eq:lambda_ER}, \eqref{eq:lambda_sum} and \eqref{eq:lambda_final} to calculate a suitable correction size \( \lambda^{i} \).

    \item Perform the adapted Euler step to obtain \(\mathbf{\bar{u}^{n+1}}\) at the next time level \( t + \Delta t\)
    \begin{equation}
        \bar{u}_{i,j}^{n+1} = \bar{u}_{i,j}^{n} - \frac{\Delta t}{h_j} \left(f_{i,j+1/2}^n - f_{i,j-1/2}^n \right) + \Delta t \lambda^{i}v_{i,j},
    \end{equation}
    where
    \begin{equation}
        \mathbf{v_{i}} = \mathbf{G}\mathbf{\bar{u}_{i}^{n}}.
    \end{equation}
\end{enumerate}
We refer to this routine by \( \mathbf{\bar{u}^{n+1}} =  \) \textbf{euler adapted} (\( \mathbf{\bar{u}^{n}}\) ) and use it to state the complete process of the modified SV method.

\textbf{Adapted SV scheme:}
\begin{enumerate}
    \item Construct the grid using Gauss-Lobatto points.
    
    \item In dependence of \(\max\limits_{{j=1,..k}}\left(h_j\right)\) and taking the CFL condition into account, determine the magnitude of one time step \(\Delta t\).
    
    \item Calculate the reconstruction matrix \( \mathbf{C} \).
    
    \item Compute the filter generator \( \mathbf{H} = \mathbf{G} \).
    
    \item Use Gauss quadrature to convert the given initial condition into a cell-averaged initial vector \( \mathbf{\bar{u}^{\left(0\right)}}\), where
    \begin{equation}
        \bar{u}_{i,j}^{\left(0\right)} = \frac{1}{h_j} \int_{x_{i,j-1/2}}^{x_{i,j+1/2}} u_0\left(x\right) dx.
    \end{equation}
    
    \item Choose a point in time \(t\) at which the numerical solution should be computed. The number of time steps \(T\) is then determined by \( \lfloor t/ \Delta t \rfloor\).

    \item For \(n = 1,...,T\) use the third order Runge-Kutta method \eqref{eq:Runge} for time integration.\\
    \begin{equation}
        \begin{aligned}
            \mathbf{\bar{u}^{\left(1\right)}} &= \textbf{euler adapted}\left( \mathbf{\bar{u}^{\left(0\right)}} \right) \\
            \mathbf{\bar{u}^{\left(2\right)}} &= \frac{3}{4}  \mathbf{\bar{u}^{\left(0\right)}} + \frac{1}{4} \cdot \textbf{euler adapted}\left( \mathbf{\bar{u}^{\left(1\right)}} \right) \\
            \mathbf{\bar{u}^{\left(3\right)}} &= \frac{1}{3} \mathbf{\bar{u}^{\left(0\right)}} + \frac{2}{3} \cdot \textbf{euler adapted}\left( \mathbf{\bar{u}^{\left(2\right)}} \right)\\
        \end{aligned}
    \end{equation}
\end{enumerate}
\clearpage
\section{Numerical Results}\label{sec:Numerical Results}

In this chapter, numerical results will be discussed to verify the adapted SV scheme presented in chapter \ref{sec:Stabilized SV scheme}. The initial section revisits the test case from chapter \ref{sec:First numerical test}, providing a comparison between the adapted scheme and the original SV scheme. Subsequently, in section \ref{sec:Numerical Results for Burgers equation}, we examine the Burgers' equation to validate the schemes performance in handling shock formation and the developement of rarefaction waves. Moving forward, section \ref{sec:Shock Tube Tests} provides an analysis of two different shock tube tests and the chapter ends with a brief convergence analysis.

Throughout this discussion, let \(N_{CV}\) and \(N_{SV}\) represent the number of control volumes and the number of spectral volumes, respectively. For all numerical tests, we will use \(N_{CV} = 4\).

\subsection{Linear Advection equation}\label{sec:Numerical Results for Linear Advection}

\begin{figure}[!h]
\begin{subfigure}{\textwidth}
    \centering
    \includegraphics[width=0.64\linewidth]{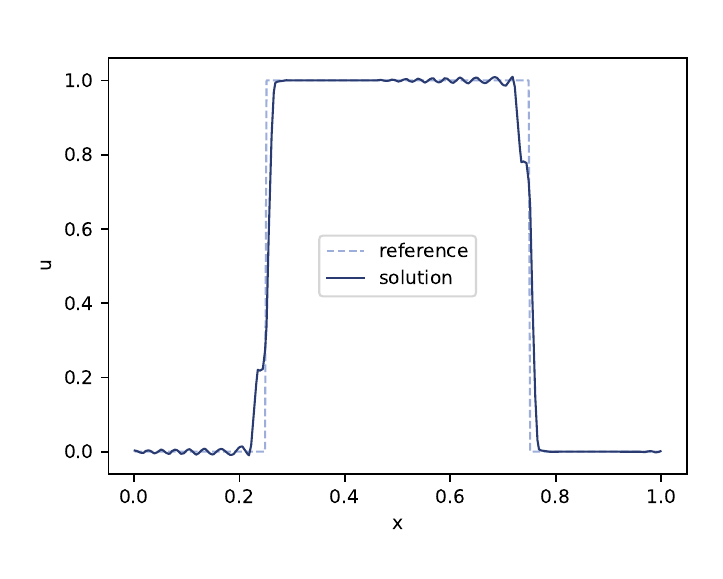}
    \caption{Adapted SV scheme: Solution at \(t = 1.0\), \(N_{SV} = 60, \, N_{CV} = 4\).}
    \label{fig:rectangle_60}
\end{subfigure}
\begin{subfigure}{\textwidth}
    \centering
    \includegraphics[width=0.64\linewidth]{rectangleunstable_CV4_SV60.pdf}
    \caption{Pure SV scheme: Solution at \(t = 1.0\), \(N_{SV} = 60, \, N_{CV} = 4\).}
    \label{fig:rectangle_unstable_comparison}
\end{subfigure}
\caption{This first numerical test demonstrates the effectiveness of adapting the SV scheme. Comparing Figures \ref{fig:rectangle_60} and \ref{fig:rectangle_unstable_comparison} shows that the adapted scheme reduces spurious oscillations while preserving sharp discontinuities. Aside from minor oscillations near the discontinuities, the solution closely matches the exact reference. In other words, the adaptation dampens unwanted oscillations without causing undesirable smoothing. This performance is likely due to the adjustable correction size, which dissipates just enough entropy without overdoing it. Overall, these results indicate a clear improvement of the adapted SV scheme over the pure SV scheme—a finding that further tests will confirm.}
\end{figure}

For the sake of comparability, our first numerical test mirrors the one in chapter \ref{sec:First numerical test}. We revisit the linear advection equation with  \( v = 1\) and a rectangle as initial condition on the domain \( \left[0, 1 \right] \):

\begin{equation}
    \partial_t u \left(x, t \right) + \partial_x u \left( x, t \right) = 0, \quad
    u_{0}\left(x\right) = 
    \begin{cases}
        1, \quad \frac{1}{4} \leq x \leq \frac{3}{4},\\
        0, \quad \text{otherwise}.
    \end{cases}
\end{equation}
Aditionally, we impose periodic boundary conditions. To apply our adapted SV scheme, a suitable entropy-entropy flux pair has to be selected. As we know from section \ref{sec:Introduction}, any convex function \( U \) can be chosen. The most straightforward and thus selected entropy is:

\begin{equation}
    U\left(u \right) = \frac{1}{2}u^{2}.
\end{equation}
The corresponding entropy flux is then given by:

\begin{equation}
    F\left( u \right) = \int U'\left( u \right) f'\left(u \right) du = \int u \, du = \frac{1}{2}u^{2}.
\end{equation}
The numerical results of the adapted SV scheme are depicted in figure \ref{fig:rectangle_60} alongside an exact reference solution at \(t = 1.0\). Again, this means that the rectangle has passed one complete round through the domain.

It is evident that the numerical solution of the stabilized SV scheme exhibits significantly fewer spurious oscillations compared to the numerical solution of the pure SV scheme presented in chapter \ref{sec:First numerical test} and depicted for comparison in figure \ref{fig:rectangle_unstable_comparison}.

\clearpage

\subsection{Burgers' equation}\label{sec:Numerical Results for Burgers equation}
Next, we will discuss numerical tests for the non-linear Burgers' equation to assess the schemes capability in handling shock formation and the developement of rarefaction waves. Burger' equation is given by
\begin{equation}
    \partial_t u\left(x, t\right) + \frac{1}{2} \partial_x u^2 \left(x, t \right) = 0
\end{equation}
For simplicity, we once again choose the entropy function

\begin{equation}
    U\left(u \right) = \frac{1}{2}u^{2}.
\end{equation}
In the case of the Burgers' equation with flux function \(f\left( u \right) = \frac{1}{2} u^{2}\), the entropy flux function is given by
\begin{equation}
    F\left( u \right) = \int U'\left( u \right) f'\left(u \right) du = \int u^{2} \, du = \frac{1}{3}u^{3}.
\end{equation}
For the first test case, we consider the problem

\begin{equation}
    \partial_t u  + \frac{1}{2} \partial_x u^{2} = 0, \quad u_{0}\left(x\right) = \sin \left( \pi x \right), \quad x \in \left[0, 2 \right]
\end{equation}
with periodic boundary conditions. From the method of characteristics, we know that a shock will be formed before \(t = 0.5\).

\begin{figure}[!h]
    \centering
    \includegraphics[width=0.65\linewidth]{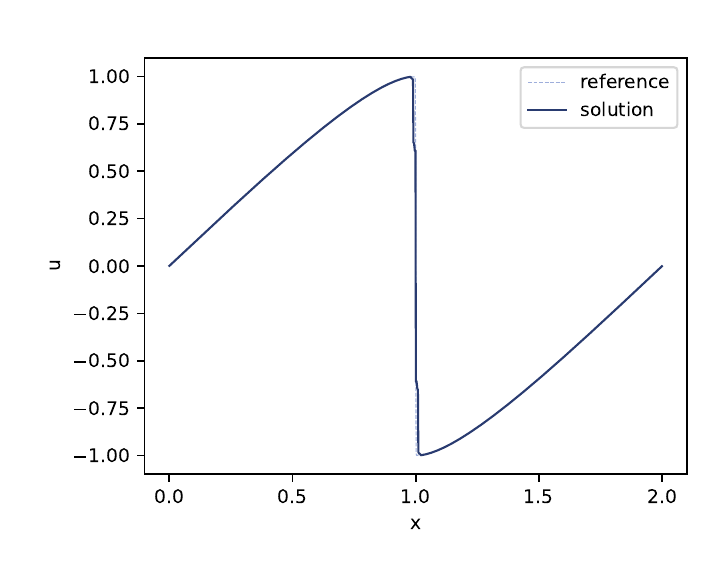}
    \caption{Solution at \(t = 0.5\), \(N_{SV} = 200, \, N_{CV} = 4\)}
    \label{fig:burger_sine_200}
\end{figure}

The numerical solution is depicted for this point in time in figure \ref{fig:burger_sine_200}. A reference solution is calculated by a Lax-Friedrichs scheme on \(3\cdot 10^{4}\) cells. The result is satisfactory as the developed shock is sharp and no spurious oscillations occur, demonstrating the schemes ability to handle shock formation.

A typical initial condition for the time evolution of a rarefaction wave is given by
\begin{equation}
    \partial_t u  + \frac{1}{2} \partial_x u^{2} = 0, \quad
    u_{0}\left(x\right) = 
    \begin{cases}
        -1, \quad  x \leq 1.0,\\
        1, \quad \text{otherwise},
    \end{cases}
\end{equation}
on the domain \( \left[0, 2 \right] \) with fixed boundary conditions. From the methods of characteristics, the analytical solution at time \(t = 1\) is known to be

\begin{equation}
    u\left(x\right) = 
    \begin{cases}
        -1, \quad  x \leq 0.5,\\
        2x - 2 , \quad  0.5 \leq x \leq 1.5,\\
        1, \quad \text{otherwise}.
    \end{cases}
\end{equation}

\begin{figure}[!h]
    \centering
    \includegraphics[width=0.65\linewidth]{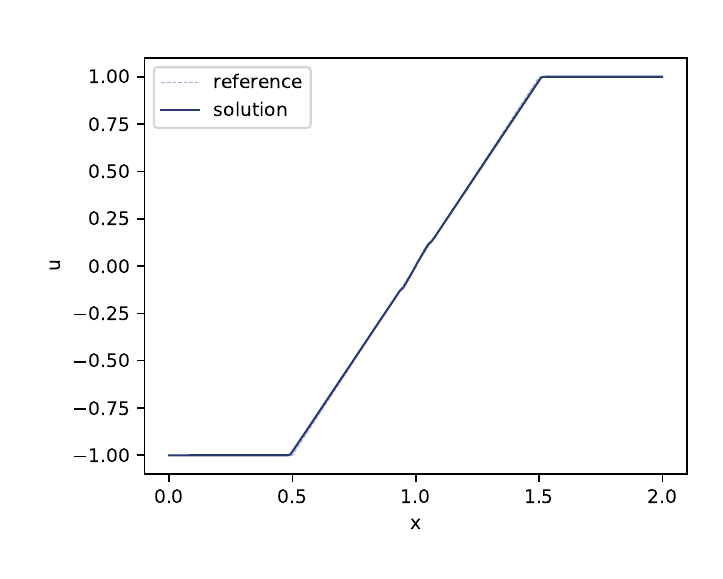}
    \caption{Solution at \(t = 0.5\), \(N_{SV} = 200, \, N_{CV} = 4\)}
    \label{fig:burger_step_200}
\end{figure}

It is shown in figure \ref{fig:burger_step_200} together with the numerical solution. The numerical result is once again satisfying as the developed rarefaction wave has no spurious smearings on the edges, indicating the schemes ability to handle the developement of rarefaction waves as well.

\clearpage

\subsection{Shock Tube Tests}\label{sec:Shock Tube Tests}
\begin{figure}[!h]
    \centering
    \includegraphics[width=0.7\linewidth]{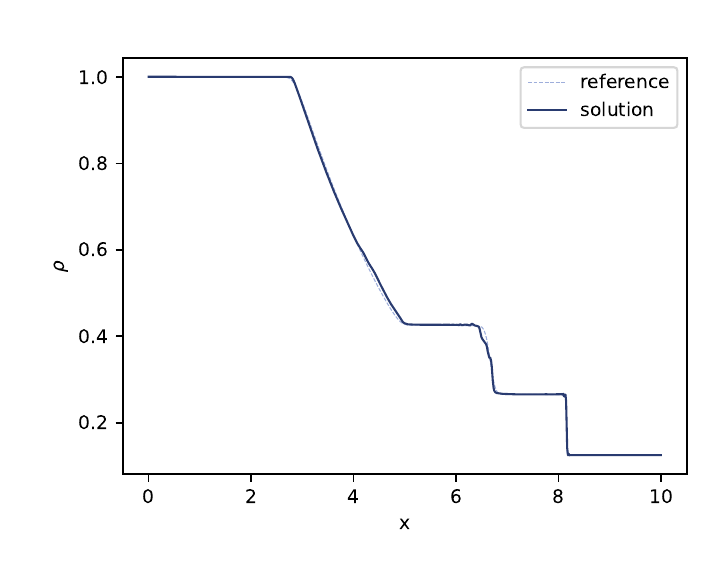}
    \caption{Density, \(N_{SV} = 200, \, N_{CV} = 4\).}
    \label{fig:st_density}
\end{figure}

\begin{figure}[!h]
    \centering
    \includegraphics[width=0.7\linewidth]{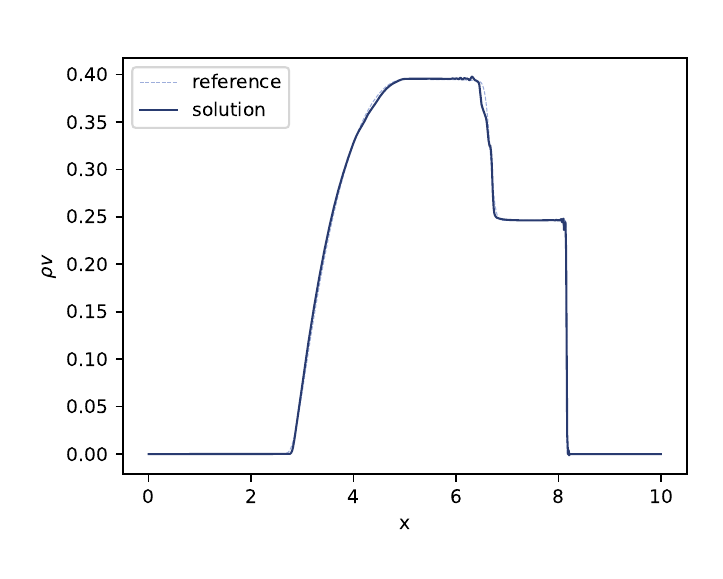}
    \caption{Momentum density, \(N_{SV} = 200, \, N_{CV} = 4\).}
    \label{fig:st_momentumdensity}
\end{figure}

\begin{figure}[!h]
    \centering
    \includegraphics[width=0.7\linewidth]{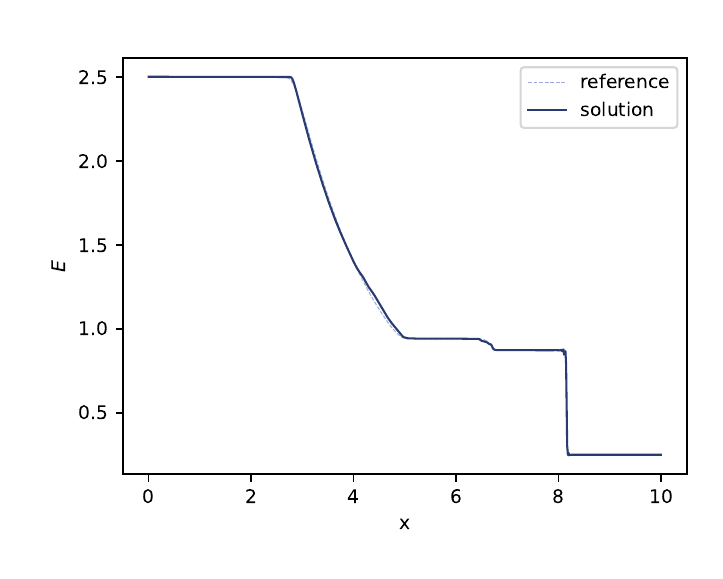}
    \caption{Energy density, \(N_{SV} = 200, \, N_{CV} = 4\).}
    \label{fig:st_energydensity}
\end{figure}

In this section, we will present numerical results for two different shock tube tests. We consider the Euler equations of gas dynamics \cite{godlewski2013numerical} with fixed boundary conditions on the domain \( \left[0, 10 \right] \):

\begin{equation}\label{eq:euler flux}
    u = \left[ \begin{array}{c} \rho\\ \rho v\\ E \end{array} \right], \quad 
    f \left(u \right) = \left[ \begin{array}{c} \rho v\\ \rho v^{2} + p\\ v \left(E + p \right) \end{array} \right], \quad
    p = \left( \gamma - 1 \right) \left(E - \frac{1}{2} \rho v^{2} \right).
\end{equation}
A siutable physical entropy-entropy flux pair is given by \cite{harten1983symmetric, tadmor2003entropy}:

\begin{equation}
    U\left(\rho, \rho v, E \right) = - \rho S, \quad F\left(\rho, \rho v, E \right) = - \rho v S, \quad S = \ln \left(p \rho ^{- \gamma} \right).
\end{equation}
In our first test case, the gas is initially at rest,  (\(v = 0\)). Additionally, we start with a higher density and pressure on the left side of \(x = 5\) than on the right. The exact initial condition \cite[Problem 6a]{shu1989efficient2} used for the numerical test is given by:

\begin{equation}\label{eq:shock tube}
    \rho_{0}\left(x\right) = 
    \begin{cases}
        1.0,\\
        0.125, 
    \end{cases}
    , \quad
    v_{0}\left(x\right) = 
    \begin{cases}
        0,\\
        0,
    \end{cases}
    , \quad
    p_{0}\left(x\right) = 
    \begin{cases}
        1.0, \quad  x < 5\\
        0.1, \quad \text{otherwise}
    \end{cases}.
\end{equation}
The initial momentum density and energy density needed for the numerical computation can be calculated from (\ref{eq:shock tube}) using \eqref{eq:euler flux}. The numerical results of the adapted SV scheme are depicted in figures \ref{fig:st_density}, \ref{fig:st_momentumdensity} and \ref{fig:st_energydensity} together with a reference solution.

\begin{figure}[!h]
    \centering
    \includegraphics[width=0.7\linewidth]{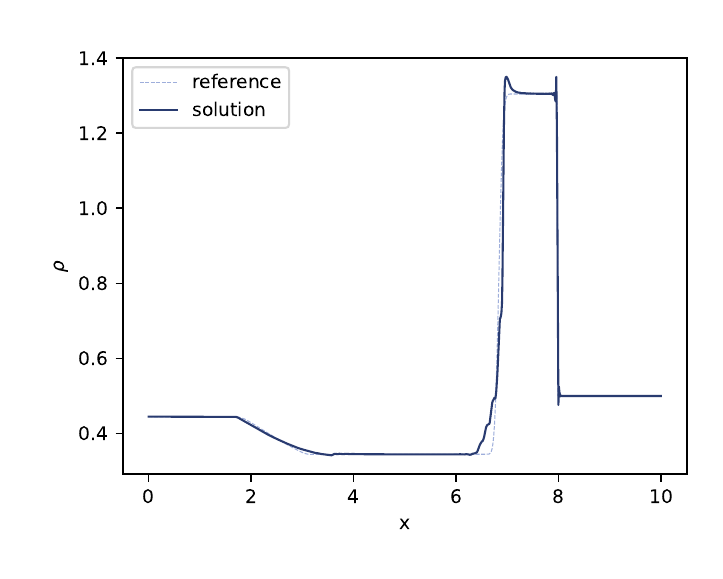}
    \caption{Density, \(N_{SV} = 200, \, N_{CV} = 4\).}
    \label{fig:st1_density}
\end{figure}

\begin{figure}[!h]
    \centering
    \includegraphics[width=0.7\linewidth]{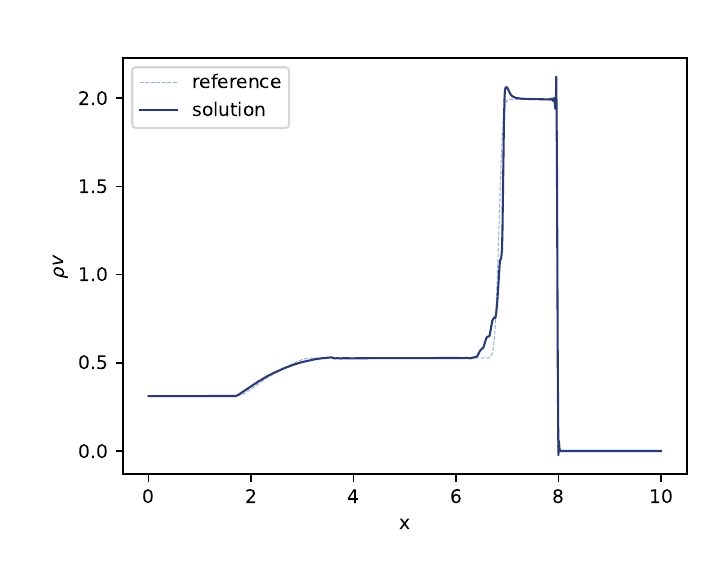}
    \caption{Momentum density, \(N_{SV} = 200, \, N_{CV} = 4\).}
    \label{fig:st1_momentumdensity}
\end{figure}

\begin{figure}[!h]
    \centering
    \includegraphics[width=0.7\linewidth]{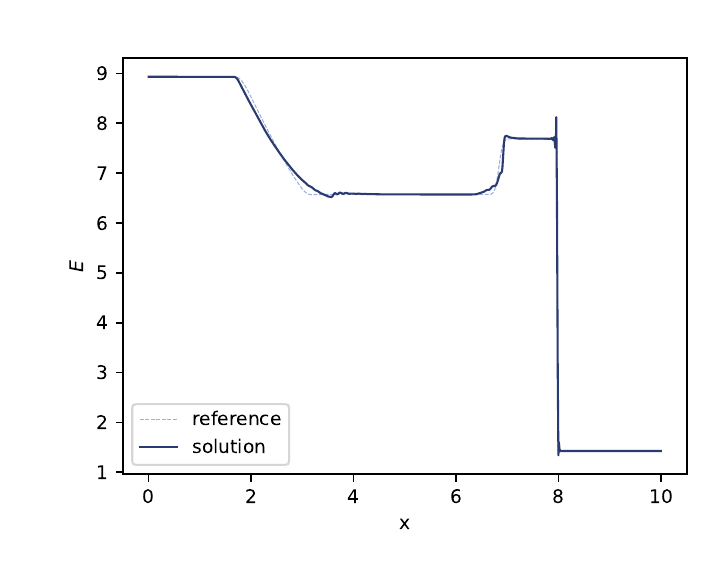}
    \caption{Energy density, \(N_{SV} = 200, \, N_{CV} = 4\).}
    \label{fig:st1_energydensity}
\end{figure}

The second shock tube test case is determined by the initial condition \cite[Problem 6b]{shu1989efficient2}:
\begin{equation}\label{eq:shock tube_1}
    \rho_{0}\left(x\right) = 
    \begin{cases}
        0.445,\\
        0.5,
    \end{cases}
    , \quad
    v_{0}\left(x\right) = 
    \begin{cases}
        0.698,\\
        0,
    \end{cases}
    , \quad
    p_{0}\left(x\right) = 
    \begin{cases}
        3.528, \quad  x < 5\\
        0.571, \quad \text{otherwise}
    \end{cases}.
\end{equation}
Note, that the initial velocity is zero only at the right side of \(x = 5\), whereas it is positive on the left side. The numerical results are shown in figures  \ref{fig:st1_density}, \ref{fig:st1_momentumdensity} and \ref{fig:st1_energydensity} alongside a reference solution.

For both shock tube tests, the reference solution was computed by a Lax-Friedrichs scheme on \(3 \cdot 10^{4}\) cells. The numerical results of the two shock tube tests look satisfactory as all shocks and contact discontinuities are sharp, and only slight overshoots and oscillations are visible directly around the discontinuities.\\

\clearpage
\subsection{Convergence Analysis}\label{sec:Convergence Analysis}
\begin{figure}[!h]
    \begin{minipage}{0.5\textwidth}
        \centering
        \includegraphics[scale = 0.95]{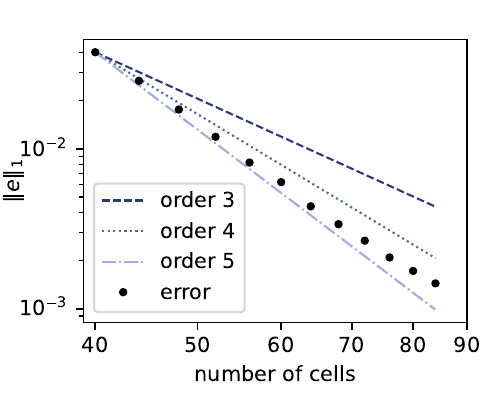}
        \subcaption{\(N_{CV} = 4\), \(L^{1}\)-norm.}
    \end{minipage}\hfill
    \begin{minipage}{0.5\textwidth}
        \centering
        \includegraphics[scale = 0.95]{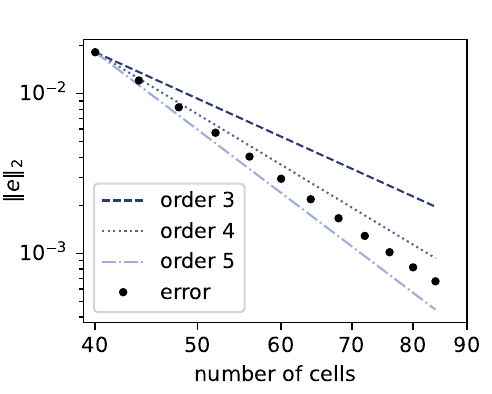}
        \subcaption{\(N_{CV} = 4\), \(L^{2}\)-norm.}
    \end{minipage}
    \caption{Convergence analysis.}
    \label{fig:convergence analysis}
\end{figure}

Finally, we will discuss the convergence speed of our modified SV scheme. As an example, we employed the modified SV scheme to calculate the numerical solution of the Euler equations with initial conditions:
\begin{equation}\label{eq:density bump}
    \rho_{0} \left( x \right) = 1 + \varepsilon \left(x \right), \quad
    v_{0} \left( x \right) = 1.0, \quad
    p_{0} \left( x \right) = 1.0, \quad
    \varepsilon \left(x \right) = e^{-\frac{\left(x - 5 \right)^{2}}{2}}.
\end{equation}
First, consider a generalized form of these initial conditions:
\begin{equation}
    \tilde{\rho_{0}} \left( x \right) = \psi \left(x \right), \quad
    \tilde{v_{0}} \left( x \right) = v_{c}, \quad
    \tilde{p_{0}} \left( x \right) = p_{c},
\end{equation}
where \(v_{c}, p_{c} \in \mathbb{R}\) and \( \psi \left(x \right)\) is some given smooth function. The solution of this initial value problem is given by
\begin{equation}
    \tilde{\rho} \left(x, t \right) = \psi\left(x-v_{c}t \right), \quad
    \tilde{v} \left( x, t \right) = v_c, \quad
    \tilde{p} \left( x, t \right) = p_c.
\end{equation}
Therefore, the special case of \(v_{c} = 1.0 \), \(p_{c} = 1.0 \) and \(\psi(x) = 1 + \varepsilon(x)\) in problem \eqref{eq:density bump} provides
\begin{equation}
    \rho \left(x, t \right) = 1 +\varepsilon \left(x - t \right), \quad
    v \left( x, t \right) = 1.0, \quad
    p \left( x, t \right) = 1.0.
\end{equation}
The numerical solution was computed at \(t = 10\) for \(N_{CV} = 4\) and \(N_{SV} = \{10, 11, 12,...,22 \}\) corresponding to \( N = \{40, 44, ..., 88\}\) cells. Subsequently, the \(L^{1}\) and \(L^{2}\) errors were calculated. They are shown in a double logarithmic plot in figure \ref{fig:convergence analysis} along with reference lines of orders 3, 4 and 5.
\vspace{-10pt}

\clearpage

Despite using a third order Runge-Kutta method for time integration, the observed convergence speed falls between orders 4 and 5, see figure \ref{fig:convergence analysis}. However, we do use a fourth-order reconstruction. This indicates a slight domination of the error in space over the error in time. Still, the observed convergence speed surpasses the expected fourth order slightly, which might be attributed to the phenomenon highlighted by Klein \cite[chapter 4.4]{klein2023stabilizing2}: When the grid is refined, the entropy dissipation estimate converges at a higher rate than the basic scheme itself. Consequently, if the correction term introduced to enforce the dissipation dominates, a higher convergence speed than expected is observed.

\clearpage

\section{Summary}\label{sec:Summary}
In this work, we modified the SV scheme proposed by Wang \cite{wang2002spectral} by incorporating an entropy stabilization method suggested by Klein for the DG method \cite{klein2023stabilizing2}. We started with an overview of the basic SV scheme and its numerical implementation, utilizing Gauss-Lobatto grid points, polynomial reconstruction, and a third-order Runge-Kutta method for time integration.\\
An initial numerical test exposed challenges with the direct applicability of the SV scheme, as spurious oscillations emerged in the numerical solution. To adress this issue, we introduced and adapted Kleins' entropy stabilization method from the DG method. Notably, the SV version of the stabilization method differed in its construction of the correction direction. As for the DG method, the direction for the SV scheme was constructed based on the time evolution of the heat equation, ensuring physical relevance. However, the SV method relies on cell averages instead of point values used in the DG version. Thus, the direction for the SV scheme was calculated by a finite volume discretization of the heat equation. This procedure appeared to be even simpler than the one used for the DG method.\\
The adapted SV scheme was then applied to various test cases to verify its practical applicability and an impovement over the basic SV scheme. Revisiting the initial test case revealed a significant reduction in spurious oscillations. Numerical tests of the Burgers' equation with different initial conditions illustrated that the adapted SV scheme can handle shock formation and the developement of rarefaction waves in the solution. Two different shock tube tests further demonstrated satisfactory results, with sharp shocks and contact discontinuities in the numerical solution and no spurious smearings at the edges of rarefaction waves. Finally, a convergence analysis revealed that the order of convergence matches the order of the original SV scheme.\\
For future enhancements, the scheme could benefit from exploring alternative spectral volume partitionings beyond the currently employed Gauss-Lobatto points. Additionally, a different discretization of the \(L^{2}\) scalar product on control volumes is worth considering.\\
Moreover, stabilization methods, different from the one introduced in this work, are undoubtedly possible to stabilize the SV scheme. For instance, stabilization methods correcting the numerical flux itself, rather than the time derivative of the state variable, could be considered \cite{vilar2019posteriori, vilar2024posteriori}.

\clearpage
\printbibliography

\end{document}